\documentclass[12pt]{amsart}
\usepackage{amsmath,amsthm}
\usepackage{amssymb}
\usepackage[all]{xypic}
\usepackage{euscript}
\begin{document} 
%
%
\swapnumbers
\newtheorem{thm}{Theorem}[section]
\newtheorem*{tha}{Theorem}
\newtheorem{lemma}[thm]{Lemma}
\newtheorem{prop}[thm]{Proposition}
\newtheorem{cor}[thm]{Corollary}
\theoremstyle{definition}
\newtheorem{claim}[thm]{Claim}
\newtheorem{defn}[thm]{Definition}
\newtheorem{example}[thm]{Example}
\newtheorem{fact}[thm]{Fact}
\theoremstyle{remark}
\newtheorem{remark}[thm]{Remark}
\newtheorem{notation}[thm]{Notation}
\newtheorem{quest}[thm]{Problem}
\numberwithin{equation}{section}
\numberwithin{figure}{section}
%
%
\newcommand{\sect}{\setcounter{thm}{0}\section}
%
%
\newcommand{\xra}[1]{\xrightarrow{#1}}
\newcommand{\xla}[1]{\xleftarrow{#1}}
\newcommand{\hra}{\hookrightarrow}
\newcommand{\adj}[2]{\substack{{#1}\\ \rightleftharpoons \\ {#2}}}
\newcommand{\hsp}{\hspace{10 mm}}
\newcommand{\hs}{\hspace{5 mm}}
\newcommand{\hsm}{\hspace{2 mm}}
\newcommand{\vs}{\vspace{7 mm}}
\newcommand{\vsm}{\vspace{2 mm}}
\newcommand{\rest}[1]{\lvert_{#1}}
\newcommand{\lra}[1]{\langle{#1}\rangle}
\newcommand{\EQUIV}{\Leftrightarrow}
\newcommand{\epic}{\to\hspace{-5 mm}\to}
\newcommand{\xepic}[1]{\xrightarrow{#1}\hspace{-5 mm}\to}
\newcommand{\hotimes}{\hat{\otimes}}
\newcommand{\hy}[2]{{#1}\text{-}{#2}}
%
%
\newcommand{\ab}{\operatorname{ab}}
\newcommand{\Aut}{\operatorname{Aut}}
\newcommand{\AutL}{\Aut_{\Lambda}}
\newcommand{\Cok}{\operatorname{Coker}\,}
\newcommand{\cof}{\operatorname{cof}}
\newcommand{\cf}{\operatorname{cf}}
\newcommand{\colim}{\operatorname{colim}}
\newcommand{\diag}{\operatorname{diag}}
\newcommand{\Ext}{\operatorname{Ext}}
\newcommand{\Fib}{\operatorname{Fib}\,}
\newcommand{\gr}{\operatorname{gr}}
\newcommand{\haut}{\operatorname{haut}}
\newcommand{\ho}{\operatorname{ho}}
\newcommand{\holim}{\operatorname{holim}}
\newcommand{\hocolim}{\operatorname{hocolim}}
\newcommand{\Hom}{\operatorname{Hom}}
\newcommand{\Id}{\operatorname{Id}}
\newcommand{\Image}{\operatorname{Im}\,}
\newcommand{\Ker}{\operatorname{Ker}\,}
\newcommand{\map}{\operatorname{map}}
\newcommand{\Obj}{\operatorname{Obj}}
\newcommand{\op}{\operatorname{op}}
\newcommand{\sk}[1]{\operatorname{sk}_{#1}}
\newcommand{\we}{\operatorname{w.e.}}
%
%
\newcommand{\A}{{\mathcal A}}
\newcommand{\hA}{\hat{\A}}
\newcommand{\Alg}[1]{{#1}\text{-}{\EuScript Alg}}
\newcommand{\Ab}{{\EuScript Ab}}
\newcommand{\Abgp}{{\Ab\Gp}}
\newcommand{\B}{{\mathcal B}}
\newcommand{\hB}{\hat{\B}}
\newcommand{\C}{{\mathcal C}}
\newcommand{\Ch}{{\EuScript Chain}}
\newcommand{\hC}{\hat{\C}}
\newcommand{\D}{{\mathcal D}}
\newcommand{\hD}{\hat{\D}}
\newcommand{\hF}{\hat{F}}
\newcommand{\G}{{\mathcal G}}
\newcommand{\Gp}{{\EuScript Gp}}
\newcommand{\Hopf}{{\EuScript Hopf}}
\newcommand{\K}{{\EuScript K}}
\newcommand{\M}{{\EuScript M}}
\newcommand{\hP}{\widehat{\pis}}
\newcommand{\Pa}{$\Pi$-algebra}
\newcommand{\PiA}{\Pi_{\A}}
\newcommand{\PihA}{\Pi_{\hA}}
\newcommand{\PiC}{\Pi_{\C}}
\newcommand{\PiD}{\Pi_{\D}}
\newcommand{\PCa}{$\PiC$-algebra}
\newcommand{\PDa}{$\PiD$-algebra}
\newcommand{\PAa}{$\PiA$-algebra}
\newcommand{\PAlg}{\Alg{\Pi}}
\newcommand{\PAAlg}{\Alg{\PiA}}
\newcommand{\PCAlg}{\Alg{\PiC}}
\newcommand{\PDAlg}{\Alg{\PiD}}
\newcommand{\hQ}{\hat{Q}}
\newcommand{\R}[1]{{\mathcal R}_{#1}}
\newcommand{\RM}[1]{{#1}\text{-}{\EuScript Mod}}
\newcommand{\Ss}{{\mathcal S}}
\newcommand{\Sa}{\Ss_{\ast}}
\newcommand{\Set}{{\EuScript Set}}
\newcommand{\Seta}{\Set_{\ast}}
\newcommand{\Tor}{\operatorname{Tor}}
\newcommand{\TT}{{\mathcal T}}
\newcommand{\TM}{{\EuScript TM}}
\newcommand{\Ta}{\TT_{\ast}}
\newcommand{\hT}{\hat{T}}
\newcommand{\hU}{\hat{U}}
%
%
\newcommand{\We}{\mathfrak{W}}
\newcommand{\fG}{\mathfrak{G}}
%
%
\newcommand{\bN}{\mathbb N}
\newcommand{\bQ}{\mathbb Q}
\newcommand{\bT}{\mathbb T}
\newcommand{\bZ}{\mathbb Z}
%
%
%
\newcommand{\bS}[1]{{\mathbf S}^{#1}}
\newcommand{\be}[1]{{\mathbf e}^{#1}}
\newcommand{\gS}[1]{{\EuScript S}^{#1}}
\newcommand{\tX}{\tilde{X}}
\newcommand{\tY}{\tilde{Y}}
%
%
\newcommand{\pis}{\pi_{\ast}}
\newcommand{\pin}{\pi^{\natural}}
\newcommand{\pinC}[1]{\pi^{\C}_{#1}}
\newcommand{\pinD}[1]{\pi^{\D}_{#1}}
\newcommand{\piA}{\pinC{\ast}}
\newcommand{\pia}[1]{\pinC{A,{#1}}}
\newcommand{\HL}[3]{H^{#1}_{\Lambda}({#2};\,{#3})}
%
%
\newcommand{\BL}{B\Lambda}
\newcommand{\BCL}{B_{\C}\Lambda}
\newcommand{\BDL}{B_{\D}\Lambda}
\newcommand{\EK}[4]{E\sp{{#1}}\sb{#2}({#3},{#4})}
\newcommand{\EL}[2]{\EK{\Lambda}{}{#1}{#2}}
\newcommand{\ECL}[2]{\EK{\Lambda}{\C}{#1}{#2}}
\newcommand{\EDL}[2]{\EK{\Lambda}{\D}{#1}{#2}}
\newcommand{\EC}[2]{\EK{}{\C}{#1}{#2}}
%
%
\newcommand{\bd}{\mathbf{d}_{0}}
\newcommand{\co}[1]{c({#1})}
\newcommand{\q}[1]{^{({#1})}}
\newcommand{\bpa}[1]{\hat{p}\q{#1}}
\newcommand{\bpc}[1]{\check{p}\q{#1}}
\newcommand{\brp}[1]{\hat{r}\q{#1}}
\newcommand{\tkp}[1]{\hat{k}_{#1}}
\newcommand{\hk}[1]{\mathbf{\hat{k}}_{#1}}
\newcommand{\Xn}[1]{X\langle{#1}\rangle}
\newcommand{\Xpn}[1]{\hat{X}\langle{#1}\rangle}
\newcommand{\hr}[1]{\hat{\rho}\q{#1}}
\newcommand{\bPa}[1]{\hat{P}_{#1}}
\setcounter{section}{-1}
%
%
\title{Comparing homotopy categories}
\author{David Blanc}
\date{June 5, 2006}
\address{Department of Mathematics, University of Haifa, 31905 Haifa, Israel}
\email{blanc@math.haifa.ac.il}
\begin{abstract}
Given a suitable functor \ $T:\C\to\D$ \ between model categories, we
define a long exact sequence relating the homotopy groups of any \
$X\in\C$ \ with those of \ $TX$, \ and use this to describe an
obstruction theory for lifting an object \ $G\in\D$ \ to $\C$.  
Examples include finding spaces with given homology or homotopy groups.
\end{abstract}

\maketitle
%
%
\section{Introduction}

A number of fundamental problems in algebraic topology can be
described as measuring the extent to which a given functor \ 
$T:\C\to\D$ \ between model categories induces an equivalence of
homotopy categories: more specifically, which objects (or maps) from
$\D$ are in the image of $T$, and in how many different ways. For example\vsm:

\begin{enumerate}
\renewcommand{\labelenumi}{\alph{enumi})\ }
\item How does one distinguish between different topological spaces
with the same homology groups, or with chain-homotopy
equivalent chain complexes?  How can one realize a given map of
chain complexes up to homotopy\vsm?
\item When do two simply-connected topological spaces have the same
rational homotopy type\vsm?  
\item When is a given topological space a suspension, up to homotopy?
Dually, how many distinct loop space structures, if any, can a given topological
space carry\vsm?
\item Is a given \Pa\ (that is, a graded group with an action
of the primary homotopy operations) realizable as the homotopy
groups of a topological space, and if so, in how many ways\vsm? 
\end{enumerate}

Our goal is to describe a unified approach to such problems that
works for functors between \emph{spherical} model categories, for which
several familiar concepts and constructions are available. These
include a set $\A$ of \emph{models} (to play the role of spheres, in
particular determining the corresponding homotopy groups \ $\piA$), \ 
Postnikov systems, and $k$-invariants. If a functor \ $T:\C\to\D$ \
respects this additional structure, we obtain a natural long exact
sequence of the form: 
\setcounter{equation}{\value{thm}}\stepcounter{subsection}
\begin{equation}\label{eles}
\dotsc\to \Gamma_{n}X~\xra{s}~\pinC{n} X~\xra{h}~
\pinD{n} TX~\xra{\partial}~\Gamma_{n-1}X\dotsc~,
\end{equation}
\setcounter{thm}{\value{equation}}
\noindent which generalizes the EHP sequence, J.H.C. Whitehead's ``certain exact
sequence'', and the spiral exact sequence of Dwyer, Kan, and
Stover. See \ \eqref{efive} \ below.

Under these hypotheses, given an object $G$ in $\D$, we want to find
an object $X$ in $\C$ with \ $TX\simeq G$. \ The key step is to choose \ 
$\piA X$ \ which fits into \ \eqref{eles}. \ We describe an inductive 
procedure for doing this, using the Postnikov systems in both
categories, together with an obstruction theory for lifting $G$ to
$\C$, along the following lines: 

%
%
\begin{tha}
Given \ $T:\C\to\D$ \ and \ $G\in\D$ \ as above, for each \ $X\in\C$ \ 
with \ $TX\simeq G$, \ there is a tower of fibrations in $\C$:
$$
\dotsb \xra{p\q{n+1}} \Xpn{n+1} \xra{p\q{n}} \Xpn{n} \xra{p\q{n-1}} \dotsb
\xra{p\q{0}} \Xpn{0}~,
$$
called the \emph{modified Postikov tower for} $X$ (Def.\ \ref{dmps}), 
with $G$ mapping compatibly to \ $T\Xpn{n}$ \ for each $n$, and \
$X\simeq\holim_{n} \Xpn{n}$.
   
Conversely, given such a tower up to level $n$, the obstruction to
extending it to level \ $n+1$ \ lies in \
$\HL{n+3}{G}{\Gamma_{n+1}\Xpn{n}}$, \ and the choices for \ $\Xpn{n+1}$ \ 
are classified by:
\begin{enumerate}
\renewcommand{\labelenumi}{$\bullet$~}
\item a class in \ $\HL{n+2}{G}{\Gamma_{n+1}{\Xpn{n}}}$;  
\item a class in \ $\HL{n+2}{\Xpn{n}}{K_{n+1}}$, \ where \ 
$K_{n+1}:=\Cok\pi_{n+2}\rho\q{n}$, \ for \ 
$\rho\q{n}:P_{n+2}G\to P_{n+2}T\Xpn{n}$.
\end{enumerate}
\end{tha}

See Theorem \ref{tfour}.

\subsection{Related work}
\label{srw}\stepcounter{thm}

The comparison problems discussed above are familiar ones in 
algebraic topology:

\begin{enumerate}
\renewcommand{\labelenumi}{\alph{enumi})\ }
\item The question of the realizability of a graded algebra as a
  cohomology ring was first raised explicitly by Steenrod in
  \cite{SteCA}, but it goes back to Hopf (in \cite{HopfT}) in the
  rational case. The ``Steenrod problem'' of realizing a given \
  $\pi_{1}$-action in homology has been studied, for example, 
  in \cite{ThP,JRSmiT1}. 
\item The comparison between integral and rational homotopy type 
  was implicit in the notion of a Serre class (cf.\ \cite{SerG,ACuH}),
  although an explicit formulation was only possible after the 
  construction of the rationalization functors of Quillen and Sullivan
  in \cite{QuiR,SulG}.
\item Possible loop space structures on a given $H$-space were
  analyzed extensively, starting with the work of Sugawara and
  Stasheff (cf.\ \cite{SugG,StaH}). The dual question on
  identifying suspensions has also been studied (see, e.g., \cite{BHilS}).
\item The question of the realizability of homotopy groups goes back
  to J.H.C.~Whitehead, in \cite{JWhR} (see also \cite{JWhSH}), and has
  reappeared in recent years in the context of \Pa s (cf.\
  \cite{DKStE,DKStB}). The relationship between homology and homotopy
  groups, which is relevant to the realization problem for both, was
  studied in  \cite{JWhSB,JWhC} (in which the ``certain exact sequence'' 
  was introduced)\vsm.
\end{enumerate}

In \cite{BauCF}, H.-J.\ Baues gave what appears to be the first 
general theory covering a wide spectrum of such realization problems.
This was an outgrowth of his earlier work on classifying
homotopy types of finite dimensional CW complexes in
\cite{BauCHF,BauHH} (which in turn builds on \cite{JWhSF}). 

His initial setting consists of a homological cofibration 
category $\C$ (corresponding to, and extending, the notion of a
resolution model category) under a theory  of coactions $\bT$
(corresponding to the category \ $\PiA$ \ of \S \ref{ssmod}). 
Baues then constructs a generalized ``certain exact sequence'' 
similar to \ \eqref{eles}, \ and
provides an inductive obstruction theory for realizing a chain complex
(or a chain map) by a $\bT$-complex (corresponding to a CW complex, or
more generally a cofibrant object in $\C$) \ -- \ see 
\cite[VI, (2.2-2.3)]{BauCF}). 

These results apply inter alia to the problem of realizing 
a chain complex by a topological space (the motivating example 
for Baues's approach), as well as to the realization of a \Pa\ 
(cf.\ \cite[D, (7.9)]{BauCF}).  However, here we consider 
functors between two different model categories that 
are not covered by \cite{BauCF}. In particular, our original
motivating example \ -- \ the realization of a \emph{simplicial} \Pa\
(by a simplicial space) \ -- \ shows that in the relative context a
more refined obstruction theory may be necessary:  
compare Theorem (2.3) of \cite[VI]{BauCF} with Theorem \ref{tfour} below.

\begin{remark}\stepcounter{subsection}
Another set of closely related questions  \ -- \ which do not quite fit
into the framework described here, though they can also be stated as
realization problems \ -- \  arise in categories of structured ring
spectra; see for example \cite{RobO} and \cite[Cor.\ 5.9]{GHopkM}\vsm . 
\end{remark}

\subsection{Notation and conventions}
\label{snac}\stepcounter{thm}

$\Ta$ \ denotes the category of pointed connected topological
spaces; \ $\Seta$ \ that of pointed sets, and \ $\Gp$ \ that of
groups. \ For any category $\C$, \ $\gr\C$ \ denotes the category of 
non-negatively graded objects over $\C$, and \ $s\C$ \ the category of
simplicial objects over $\C$. \ $s\Set$ \ is denoted by $\Ss$, \
$s\Seta$ \ by \ $\Sa$, \ and \ $s\Gp$ \ by \ $\G$. \ The constant
simplicial object an an object \ $X\in\C$ \ is written \ $\co{X}\in s\C$. 

If $\C$ has all coproducts, then given \ $A\in\Ss$ \ and \ $X\in\C$, \
we define \ $X\hotimes A\in s\C$ \ by \ 
$(X\hotimes A)_{n}:=\coprod_{a\in A_{n}} X$, \ 
with face and degeneracy maps induced from those of $A$. 
For \ $Y\in s\C$, \ define \ $Y\otimes A\in s\C$ \ by \ 
$(Y\otimes A)_{n}:= \coprod_{a\in A_{n}} Y_{n}$ \ (the diagonal of the
bisimplicial object \ $Y\hotimes A$) \ -- \ so that for \ $X\in\C$ \
we have \ $X\hotimes A=\co{X}\otimes A$.

The category of chain complexes of $R$-modules is denoted by \
$\Ch_{R}$ \ (or simply \ $\Ch$, \ for $R=\bZ$).

\subsection{Organization:}
\label{sorg}\stepcounter{thm}

In Section \ref{csmc} we define \emph{spherical} model categories,
having the additional structure mentioned above. 
Most examples of such categories are in particular \emph{resolution} 
model categories, which are described in Section \ref{crmc}; we
explain how to produce the needed structure for them in Section \ref{csp}.
We define \emph{spherical functors} between such categories, and
construct the comparison exact sequence for them, in
Section \ref{csf}. This is applied in Section \ref{ccps} to study the
effect of a spherical functor on Postnikov systems. Finally, in
Section \ref{cfib} we construct an obstruction theory as above for the
fiber of a spherical functor. In Section \ref{cat} we indicate how the
theory works for the above examples.

\subsection{Acknowledgements}
\label{sack}\stepcounter{thm}

I would like to thank Paul Goerss for many hours of discussion on various
issues connected with this paper, and especially for his essential
help with Sections \ref{ccps}-\ref{cfib}, the technical core of this
note. I would also like to thank Hans Baues for explaining the
relevance of his work in \cite{BauCF} to me.

%
%
\sect{Spherical model categories}
\label{csmc}

Before defining the additional structure we shall need, we briefly
recapitulate the relevant homotopical algebra:

\subsection{Model categories}
\label{smc}\stepcounter{thm}

Recall that a \emph{model category} is a bicomplete category $\C$
equipped with three classes of maps: weak equivalences,
fibrations, and cofibrations, related by appropriate lifting
properties. By inverting the weak equivalences we obtain the
associated homotopy category \ $\ho\C$, \ with morphism set \ 
$[X,Y]=[X,Y]_{\C}$. \ We shall concentrate on \emph{pointed} model
categories (with null object $\ast$). \ See \cite{QuiH} or \cite{PHirM}.   

\subsection{The set of models}
\label{ssmod}\stepcounter{thm}

The additional initial data that we shall require for our model
category consists of a set $\A$ of cofibrant homotopy cogroup objects
in $\C$, called \emph{models} (playing the role of the spheres in $\Ta$).
Given such a set $\A$, let \ $\PiA$ \ denote the smallest subcategory of $\C$
containing $\A$ and closed under weak equivalences, arbitrary
coproducts, and suspensions. Note that every object in \ $\PiA$ \ is a
homotopy cogroup object, too.  

\begin{example}\label{esmod}\stepcounter{subsection}
Let \ $\C=\G$ \ be the category of simplicial groups, \
$S^{k}=\Delta[k]/\partial\Delta[k]$ \ the standard simplicial $k$-sphere
in\ $\Sa$, \ $G:\Sa\to\G$ \ the Kan's loop functor (cf.\ \cite[\S 26.3]{MayS}), 
and \ $F:\Sa\to\G$ \ the free group functor. For each \ $n\geq 1$,  \  
$\gS{n}:= GS^{n}\in\G\cong FS^{n-1}$ \ will be called the
\textit{$n$-dimensional $\G$-sphere}, \ with \
$\Sigma^{k}\gS{n}\simeq\gS{n+k}$. \  These, and their coproducts, are
cofibrant strict cogroup objects for $\G$. \ Here \
$\A:=\{\gS{1}=\co{\bZ}\}$; \ in fact, throughout this paper $\A$ will
be either a singleton, or countable. 
\end{example}

\begin{remark}\label{ssg}\stepcounter{subsection}
The adjoint pairs of functors:
$$
\Ta\ \ \substack{S\\ \rightleftharpoons\\ \|-\|}\ \ \Sa\ \ 
\substack{G\\ \rightleftharpoons\\ \bar{W}}\ \ \G
$$
\noindent induce equivalences of the corresponding homotopy categories \ -- \ 
where \ $\bar{W}:\G\to\Sa$ \ is the 
Eilenberg-Mac~Lane classifying space functor, \ $S:\Ta\to\Sa$ \ is the 
singular set functor, and \ $\|-\|:\Sa\to\Ta$ \ is the geometric realization 
functor (cf.\ \cite[\S 14,23]{MayS}). Thus to study the usual homotopy category
of (pointed connected) topological spaces, we can work in $\G$ (or \ $\Sa$), \ 
rather than \ $\Ta$. 
\end{remark}

\begin{defn}\label{dpis}\stepcounter{subsection}
If $\A$ is a set of models for $\C$, then given \ $X\in\C$, \  for
each \ $A\in\A$ \ let \ $\pia{k}(X):=[\Sigma^{k}A,X']_{\C}$, \
where \ $X'\to X$ \ is a (functorial) fibrant replacement. We write \ 
$\pinC{k} X$ \ for \ $(\pia{k}X)_{A\in\A}$, \ and \ 
$\piA X:=(\pinC{k}X)_{k=0}^{\infty}$.
\end{defn}

\subsection{Theories and algebras}
\label{staa}\stepcounter{thm}

Recall that a \emph{theory} is a small category  $\Theta$ with 
finite products (so in particular, an FP-sketch \ -- \ cf.\ 
\cite[\S 5.6]{BorcH2}), and a $\Theta$-\emph{algebra} (or
\emph{model}) is a product-preserving functor \ $\Theta\to\Set$. \ 
Think of $\Theta$ as encoding the operations and relations for 
a ``variety of universal algebras'', the category \ $\Alg{\Theta}$ \ 
of $\Theta$-algebras (which is \emph{sketched} by $\Theta$).

For example, the obvious category $\fG$, which sketches groups, is
equivalent to the opposite of the homotopy category of (finite) wedges
of circles. An $\fG$-\emph{theory} $\Theta$
(cf.\ \cite[\S 2]{BPescF}) is one equipped with a map of theories \ 
$\coprod_{S}\,\fG\to\Theta$ \ (coproduct taken in the category of
theories, over some index set $S$) which is bijective on objects. This
implies that each $\Theta$-algebra has the underlying 
structure of an $S$-graded group, so that \ $\Alg{\Theta}$ \  
can be thought of as a ``variety of (graded) groups with operators'' 
(cf.\ \cite[I, (2.5)]{BauCF}).

\begin{remark}\label{rpis}\stepcounter{subsection}
We will assume that all the functors \ $\pinC{n}$ \ ($n\geq 0$) \ take
value in a category \ $\PCAlg$ \ sketched by a $\fG$-theory $\Theta$, 
and thus equipped with a faithful forgetful
functor \ $U_{\C}:\PCAlg\to\Gp^{\A}$ \ into the category of
$\A$-graded groups. 
The objects of \ $\PCAlg$ \ are called \ \emph{\PCa s}.

For topological spaces, with \ $\A=\{\bS{1}\}$, \ the \PCa s are simply
groups. If we use rational spheres as the models, then \ $\PCAlg$ \ 
is the category of $\bQ$-vector spaces. A more interesting example
appears in \S \ref{dpa} below. 
\end{remark}

\subsection{Constructions based on models}
\label{scbm}\stepcounter{thm}

There are a number of familiar constructions for topological spaces which we
require for our purposes.  We can \emph{define} them once we are given
a set of models $\A$ as above, although they do not always exist (see
\S \ref{snsrmc} below).  

\begin{defn}\label{dfpt}\stepcounter{subsection}
A \emph{Postnikov tower} (with respect to $\A$) is a functor that
assigns to each \ $Y\in\C$ \ a tower of fibrations: 
$$
\dotsc \to P^{\A}_{n}Y\xra{p\q{n}}P^{\A}_{n-1}Y\xra{p\q{n-1}}\dots 
\to P^{\A}_{0}Y~,
$$
as well as a weak equivalence \ $r:Y\to P^{\A}_{\infty}Y:=\lim_{n}P^{\A}_{n}Y$ \ 
and fibrations \ $P^{\A}_{\infty}Y\xra{r\q{n}}P^{\A}_{n}Y$ \ such that \
$r\q{n-1}=p\q{n}\circ r\q{n}$ \ for all $n$. \ Finally, \ 
$(r\q{n}\circ r)_{\#}:\pinC{k}Y\to\pinC{k}(P^{\A}_{n}Y)$ \ 
is an isomorphism for \ $k\leq n$, \ and \ $\pinC{k}(P^{\A}_{n}Y)$ \
is zero for \ $k>n$. \ 

When $\A$ is clear from the context, we denote \ $P_{n}^{\A}$ \
simply by \ $P_{n}$.
\end{defn}

\begin{example}\label{epost}\stepcounter{subsection}
For a free chain complex \ $C_{\ast}\in\Ch_{R}$ \ of modules over a
ring $R$, we may take \ $C'_{\ast}:=P_{n}C_{\ast}$ \ where \
$C'_{i}=C_{i}$ \ for \ $i\leq n+1$, \ $C'_{n+2}=Z_{n+1}C_{\ast})$, \ and \ 
$C'_{i}=0$ \ for \ $i\geq n+3$. \ The map \ $r\q{n}:C_{\ast}\to C'_{\ast}$ \ 
is defined by \ $r\q{n}_{n+2}:=\partial_{n+2}:C_{n+2}\to Z_{n+1}C_{\ast}$. 
\end{example}

\begin{defn}\label{drem}\stepcounter{subsection}
Given an \PCa\ \ $\Lambda$, \ a \emph{classifying object} \ $\BCL$ \
(or simply \ $\BL$) \ for $\Lambda$ is any \ $B\in s\C$ \ such that \
$B\simeq P_{0}K$ \ and \ $\pinC{0}B\cong\Lambda$.
\end{defn}

The name is used by analogy with the classifying space of a group,
which classifies $G$-bundles. One can interpret \ $\BCL$ \ similarly,
though perhaps less naturally (see, e.g., \cite[\S 4.6]{BJTurR}).

\begin{defn}\label{dmod}\stepcounter{subsection}
A  \emph{module over a \PCa} $\Lambda$ is an abelian group object 
in \ $\PCAlg/\Lambda$ (cf.\ \cite[\S 2]{QuiC}), and the
category of such is denoted by \ $\RM{\Lambda}$. 
\end{defn}

\begin{remark}\label{rmodule}\stepcounter{subsection}
Since any \PCa\ is in particular a (graded) group, if \
$p:Y\to\Lambda$ \ is a module, then \ $Y=K\times\Lambda$  \ 
(as sets!) for \ $K:=\Ker(p)$, \ with an appropriate \PCa\ structure 
(cf.\ \cite[\S 3]{BlaG}). \ We may call $K$ itself a
$\Lambda$-\emph{module} (which corresponds to the traditional
description of an $R$-module, for a ring $R$). 
\end{remark}

\begin{example}\label{emodule}\stepcounter{subsection}
For any object \ $X\in\C$ \ as above, the \ $\A\times\bN$-graded group \ 
$\piA X$ \ has an action of the $\A$-primary homotopy
operations, corepresented by the maps in \ $\ho\PiA$ \ (see 
\S \ref{dpa} below). In particular, one of these operations,
corresponding to the action of the fundamental group on the higher
homotopy groups, makes each \ $\pinC{n}X$ \ ($n\geq 1$) \ into a
module over $\pinC{0}X$ \ (see Fact \ref{ffour} below).
\end{example}

\begin{defn}\label{dem}\stepcounter{subsection}
Given an abelian \PCa\ $M$ and an integer \ 
$n\geq 1$, \ an \emph{$n$-dimensional $M$-Eilenberg-Mac~Lane object} \ 
$\EC{M}{n}$ \ (or simply \ $\EK{}{}{M}{n}$) \ is any \ $E\in s\C$ \ such that \ 
$\pinC{n}E\cong M$ \ and \ $\pinC{k}E=0$ \ for \ $k\neq n$.
\end{defn}

\begin{defn}\label{deaem}\stepcounter{subsection}
Given a \PCa\ $\Lambda$, \ a module $M$ over $\Lambda$, and an integer \ 
$n\geq 1$, \ an \emph{$n$-dimensional extended $M$-Eilenberg-Mac~Lane object} \ 
$\ECL{M}{n}$ \ (or simply \ $\EL{M}{n}$) \ is any homotopy abelian group object \ 
$E\in s\C/\Lambda$, \ equipped with a section $s$ for \ 
$p\q{0}:E\to P_{0}E\simeq\BL$, \ such that \ 
$\pinC{n}E\cong M$ \ as modules over $\Lambda$; \ and \ $\pinC{k}E=0$ \ 
for \ $k\neq 0,n$.
\end{defn}

\begin{defn}\label{dki}\stepcounter{subsection}
Given a Postnikov tower functor as in \S \ref{dfpt}, an $n$-th 
\emph{$k$-invariant square} (with respect to $\A$)  
is a functor that assigns to each \ $Y\in\C$ \ a homotopy pull-back square:
\setcounter{equation}{\value{thm}}\stepcounter{subsection}
\begin{equation}\label{ezero}
\xymatrix@R=25pt{\ar @{} [dr] |<<<{\framebox{\scriptsize{PB}}}
P^{\A}_{n+1}Y \ar[r]^{p\q{n+1}} \ar[d] &
  P^{\A}_{n}Y \ar[d]^{k_{n}}\\ \BL \ar[r] & \EL{M}{n+2}
}
\end{equation}
\setcounter{thm}{\value{equation}}
\noindent for \ $\Lambda:=\pinC{0}Y$ \ and \ $M:=\pinC{n+1}Y$, \ where
\ $p\q{n+1}:P_{n+1}Y\to P_{n}Y$ \ is the given fibration of the
Postnikov tower. 

The map \ $k_{n}:P_{n}Y\to\EL{M}{n+2}$ \ is the $n$-th
(functorial) $k$-\emph{invariant} for $Y$. 
\end{defn}

\begin{example}\label{ekinv}\stepcounter{subsection}
If \ $C_{\ast}$ \ is a chain complex of $R$-modules, and \ 
$P_{n}C_{\ast}=C'_{\ast}$ \ as in \S \ref{epost}, we may take \
$\EK{}{}{H_{n+1}C_{\ast}}{n+2}=E_{\ast}$, \ where \ $E_{i}=0$ \ for \
$i<n+2$, \ $E_{n+2}=Z_{n+1}C_{\ast}$, \ and \  
$E_{n+3}=B_{n+1}C_{\ast}$. \ Then \ $k_{n}:C'_{\ast}\to E_{\ast}$ \ is
defined by \ $\Id:C'_{n+1}\to E_{n+1}$.

Of course, if $R$ is a principle ideal domain (or a hereditary ring),
such as $\bZ$, then the $k$-invariants for \ $C_{\ast}$ \ are trivial,
since in that case any two free (or projective) chain complexes with
the same homology are homotopy equivalent, by \cite[Prop.\ 3.5]{DolH}. 
But this need not hold for an arbitrary ring $R$.
\end{example}

\subsection{Spherical models}
\label{dsmod}\stepcounter{thm}

A set of objects \ $\A:=\{A\}_{A\in\A}$ \ in a model category
$\C$ is called a collection of \emph{spherical models} if the
following axioms hold\vsm : 

\begin{enumerate}
\renewcommand{\labelenumi}{Ax~\arabic{enumi}.~}
\item Each \ $\Sigma^{n}A$ \ ($A\in\A$, \ $n\in\bN$) \ is a
  cofibrant homotopy cogroup object in $\C$\vsm.
\item For any \ $X\in\C$ \ and \ $n\geq 1$, \ $\pinC{n}X$ \ has a
  natural structure of a module over \ $\pinC{0}X$\vsm.
\item A map \ $f:X\to Y$ is a weak equivalence if and only if \ $\pia{n}f$ \ is a 
weak equivalence for each \ $A\in\A$ \ and \ $n\in\bN$\vsm.
\item $\C$ has Postnikov towers with respect to $\A$\vsm.
\item For every \PCa\ $\Lambda$ \ and module $M$ over $\Lambda$, the
  classifying object \ $\BL$ \ and extended $M$-Eilenberg-Mac~Lane
  object \ $\EL{M}{n}$ \ exist (and are unique up to homotopy) for
  each \ $n\geq 1$\vsm.
\item $\C$ has $k$-invariant squares with respect to $\A$ for
  each \ $n\geq 0$\vsm. 
\end{enumerate}

If each model \ $\Sigma^{k}A$ \ ($A\in\A$, \ $k\in\bN$) \ is a 
cofibrant \emph{strict} cogroup object \ -- \ which implies that 
every object in \ $\PiA$ \ is such, up to weak equivalence \ -- \ we
call $\A$ a collection of \emph{strict} spherical models.

A pointed simplicial model category \ $\C$ equipped with
a collection \ $\A:=\{A\}_{A\in\A}$ \ of spherical models
is called a \emph{spherical model category}, \ and we denote it
by \ $\lra{\C;\A}$. \ Such a category is \emph{stratified} in the
sense of Spali\'{n}ski (cf.\ \cite{SpalSM}). 

\begin{example}\label{esmc}\stepcounter{subsection}
The category \ $\Sa$ \ of pointed simplicial sets, as well as the
category \ $\Ta$ \ of pointed connected topological spaces, have
spherical model category structures with \ $\A=\{\bS{1}\}$. \ 
(Functorial $k$-invariants in these categories are provided by the
construction of \cite[\S 5]{BDGoeR}; \ in both cases \
$\PCAlg\approx\Gp$).  \ Similarly for the category \ $\Ch_{R}$ \ of
chain complexes over $R$, with the constructions indicated in \S \ref{epost} and
\S \ref{ekinv}.
\end{example}

In the examples we have in mind, our model categories enjoy
additional useful properties, which we can summarize in the following:

\begin{defn}\label{dssmc}\stepcounter{subsection}
A spherical model category \ $\lra{\C;\A}$ \ as above is called
\emph{strict} if the following axioms hold\vsm :  

\begin{enumerate}
\renewcommand{\labelenumi}{Ax~\arabic{enumi}.~}
\item $\C$ \ is a pointed right-proper cofibrantly generated simplicial 
model category (cf.\ \cite[11.1, 13.1]{PHirM}), in which every object is fibrant.
\item $\C$ is equipped with a faithful forgetful functor \
$\hat{U}:\C\to\D$, \ with left adjoint $\hat{F}$ \ -- \ where \ $\D$
  is one of the ``categories of groups'' \ $\D=\Gp$, \ $gr\Gp$, \
  $\G$, \ $\RM{R}$, \ or \ $s\RM{R}$, \ for some ring $R$.
\item The adjoint pair \ $(\hat{U},\hat{F})$ \ \emph{create} the 
model category structure on $\C$ in the sense of \cite[\S 4.13]{BlaN} \ -- \  
so in particular $\hat{U}$ creates all limits in $\C$.
\item $\A$ is a collection of strict spherical models, each of which lies
in the image of the composite \ $\hat{F}\circ F':\Ss\to \C$, \ 
where \ $F':\Ss\to\D$ \ is adjoint to the forgetful functor \
$U':\D\to\Ss$, \ with the group structure on \ $\Hom_{\C}(A,X)$ \
induced from that of \ $\hat{U}(X)$.
\end{enumerate}
\end{defn}

%
%
\sect{Resolution model categories}
\label{crmc}

Many examples of spherical model categories fit into the
framework originally conceived by Dwyer, Kan and Stover in
\cite{DKStB} under the name of ``$E^{2}$ model categories,'' and later 
generalized by Bousfield (see \cite{BousCR,JardBE}. A slightly
different generalization is given by Baues in 
\cite[Ch.~D, \S 2]{BauCF} under the name of \emph{spiral} model categories.

First, some preliminary concepts:  
 
\begin{defn}\label{dmat}\stepcounter{subsection}
The $n$-th \emph{matching object} for a simplicial object $X$ over
$\C$ is defined by
$$
M_{n}X=\{(x_{0},\dotsc,x_{n})\in(X_{n-1})^{n+1}~|\ 
d_{i}x_{j}=d_{j-1}x_{i} \text{\ \ for all\ }0\leq i<j\leq n\}
$$
(see \cite[X,\S 4.5]{BKaH}). Note that each face map \
$d_{k}:X_{n}\to X_{n-1}$ \ factors through the obvious map \ 
$\delta_{n}:X_{n}\to M_{n}X$.
\end{defn}

\begin{defn}\label{dlat}\stepcounter{subsection}
The $n$-th \emph{latching object} of a simplicial object $X$ over $\C$
is defined \ $L_{n}X:=\coprod_{0\leq i\leq n-1} X_{n-1}/\sim$, \ where for any \ 
$x\in X_{n-k-1}$ \ and \ $0\leq i\leq j \leq n-1$ \ we set \ 
$s_{j_{1}}s_{j_{2}}\dotsc s_{j_{k}}x$ \ in the $i$-th copy of \ $X_{n-1}$ \ 
equivalent to \ $s_{i_{1}}s_{i_{2}}\dotsc s_{i_{k}}x$ \ in the $j$-th copy 
of \ $X_{n-1}$ \ whenever the simplicial identity \ 
$s_{i}s_{j_{1}}s_{j_{2}}\dotsc s_{j_{k}}=s_{j}s_{i_{1}}s_{i_{2}}\dotsc
s_{i_{k}}$ \ holds. The map \ $\sigma_{n}:L_{n}X\to X_{n}$ \ is defined
\ $\sigma_{n}(x)_{i}=s_{i}x$, \ where \ $(x)_{i}\in (X_{n-1})_{i}$.
\end{defn}

There are two canonical ways to extend a given model category structure on 
$\hC$ to \ $\C:=s\hC$:

\subsection{The Reedy model structure}
\label{sreed}\stepcounter{thm}

This is defined by letting a simplicial map \ $f:X\to Y$ \ in \
$\C:=s\hC$ \ be:

\begin{enumerate}
\renewcommand{\labelenumi}{(\roman{enumi})}
\item a weak equivalence if \ $f_{n}:X_{n}\to Y_{n}$ \ is a weak equivalence 
in $\hC$ \ for each \ $n\geq 0$;
\item a (trivial) cofibration if \ 
$f_{n}\amalg \sigma_{n}:X_{n}\amalg_{L_{n}X}L_{n}Y\to Y_{n}$ \ is a 
(trivial) cofibration in $\hC$ \ for each \ $n\geq 0$;
\item a (trivial) fibration if \ 
$f_{n}\times \delta_{n}:X_{n}\to Y_{n}\times_{M_{n}Y}M_{n}X$ \ is a 
(trivial) fibration in $\hC$ \ for each \ $n\geq 0$.
\end{enumerate}

See \cite[15.3]{PHirM}.

\subsection{The resolution model category}
\label{srmc}\stepcounter{thm}

Let $\hC$ be a pointed cofibrantly generated right proper model category
(in our cases, every object will be fibrant, though this is not needed
in general  \ -- \ cf.\ \cite{JardBE}). Given a 
set $\hA$ of models for $\hC$ \ (\S \ref{ssmod}), we let \ 
$\A:=\{\co{\Sigma^{k}\hat{A}}\}_{k\in\bN,\hat{A}\in\hA}$ \ 
(the constant simplicial objects on \ $\Sigma^{k}\hat{A}\in\hC$) \ be
the set of models for $\C$. \ Note that \ 
$\Sigma^{n}\co{\Sigma^{k}\hat{A}}:=\co{\Sigma^{k}\hat{A}}\hotimes S^{n}$ \ 
(\S \ref{ssmod}), \ so we shall generally reserve the notation \
$\Sigma^{k}$ \ for (internal) suspension in $\hC$, and \
$-\otimes S^{n}$ \ for the (simplicial) suspension in \ $\C=s\hC$. 

\begin{remark}\label{rmods}\stepcounter{subsection}
If we do not assume that each \ $\hat{A}\in\hA$ \ is a homotopy
cogroup object in $\hC$, we take \
$\A:=\{(\Sigma^{k}\hat{A})\otimes S^{1}\}_{A\in\A}$ \  as our
collection of models for $\C$.  
\end{remark}

\begin{defn}\label{dffm}\stepcounter{subsection}
A map \ $f:V\to Y$ \ in  \ $\C=s\hC$ \ is called \emph{homotopically
$\hA$-free} if for each \ $n\geq 0$, \ there is 

\begin{enumerate}
\renewcommand{\labelenumi}{\alph{enumi})\ }
\item a cofibrant object \ $W_{n}$ \ in \ $\Pi_{\hA}\subset\hC$, \ and
\item a map \ $\varphi_{n}:W_{n}\to Y_{n}$ \ in $\C$ inducing a trivial 
cofibration \ $(V_{n}\amalg_{L_{n}V}L_{n}Y)\amalg W_{n}\to Y_{n}$.
\end{enumerate}
\end{defn}

We define the \emph{resolution model category structure} on \ $s\hC$ \ 
determined by $\hA$, by letting a simplicial map \ $f:X\to Y$ \ be:

\begin{enumerate}
\renewcommand{\labelenumi}{(\roman{enumi})}
\item a \emph{weak equivalence} if \ $\pia{n}f$ \ is a weak equivalence 
of simplicial groups for each \ $A\in\A$ \ and \ $n\geq 0$.
\item a \emph{cofibration} if it is a retract of a homotopically $\hA$-free map;
\item a \emph{fibration} if it is a Reedy fibration (\S \ref{sreed}(iii)) 
and \ $\pia{n}f$ \ is a fibration of simplicial groups for each \ 
$A\in\A$ \ and \ $n\geq 0$
\end{enumerate}

\begin{defn}\label{dnc}\stepcounter{subsection}
Given a fibrant \ $X\in s\hC$, \ define its $n$-\emph{cycles} object \ 
$Z_{n}X$ \ to be \ $\{ x\in X_{n}\,|\ d_{i}x=0 \ \text{for}\
i=0,\dotsc,n\}$ \ (the fiber of \ $\delta_{n}:X_{n}\to M_{n}X$ \ of 
\S \ref{dmat}). \ Similarly, the $n$-\emph{chains} object for $X$ is \ 
$C_{n}X=\{ x\in X_{n}\,|\ d_{i}x=0 \ \text{for}\ i=1,\dotsc,n\}$.
\end{defn}

If $X$ is fibrant, the map \ 
$\bd=\bd^{n}:=d_{0}\rest{C_{n}X}:C_{n}X\to Z_{n-1}X$ \ fits into a 
fibration sequence:
%
\setcounter{equation}{\value{thm}}\stepcounter{subsection}
\begin{equation}\label{etwo}
\dotsb \Omega Z_{n}X\to Z_{n+1}X\xra{j^{X}_{n+1}}C_{n+1}X \xra{\bd^{n+1}} Z_{n}X
\end{equation}
\setcounter{thm}{\value{equation}}
\noindent (see \cite[Prop.\ 5.7]{DKStB}). 

\begin{defn}\label{dpa}\stepcounter{subsection}
A \emph{\PAa} \ is a product-preserving functor from \
$(\ho\PiA)^{\op}$ \ to sets. The category of \PAa s is denoted by \ $\PAAlg$.

Equivalently, we can think of an \PAa\ $\Lambda$ as an \ 
$\bN\times\A$-graded group equipped with an action of the $\A$-primary
homotopy operations (corepresented by the maps in \ $\ho\PiA$). 

Thus we can think of the functor \ $\piA$ \ as taking value in \ $\PAAlg$.
This explains the additional \PCa\ structure on the $\A$-graded groups \ 
$\pinC{n}X$, \ mentioned in \S \ref{rpis}: when \ $\C=s\hC$, \ we
have \ $\PCAlg:=\Alg{\PihA}$.
\end{defn}

\begin{example}\label{epa}\stepcounter{subsection}
When \ $\C=\G$, \ and \ $\A=\{\gS{1}\}$ \ -- \ so \ $\PiA$ \ is the
category of wedges of $\G$-spheres (\S \ref{esmod}) \ -- \  then (up 
to indexing) \ $\PAAlg$ \ is the usual category of \Pa s (see
\cite[\S 2]{StoV}): \ graded groups equipped with an action of the
primary homotopy operations (Whitehead products and compositions). 
\end{example}

\subsection{Examples of resolution model categories}
\label{ermc}\stepcounter{thm}

In this paper we shall be interested mainly in the following instances
of resolution model categories:

\begin{enumerate}
\renewcommand{\labelenumi}{(\alph{enumi})\ }
\item Let \ $\hC=\Gp$ \ with the \emph{trivial} model category structure: \ 
i.e., only isomorphisms are weak equivalences, and every map is both a fibration
and a cofibration. Let \ $\hA=\{\bZ\}$ \ consist of the free cyclic group 
(whose coproducts are the cogroup objects in \ $\Gp$). \ 
The resulting resolution model category
structure on \ $\G:= s\hC$ \ is the usual one (cf.\ \cite[II,\S 3]{QuiH}. \ 
Here \ $\PCAlg\approx \Gp$ \ -- \ there is no extra structure on the
individual homotopy groups of a simplicial group.

Note that if we tried to do the same for \ $\hC=\Set$, \ there are no 
nontrivial cogroup objects, while in \ $\Ss$ \ not all objects are fibrant.
Note also that the category \ $\Ta$ \ of pointed topological spaces, 
which is one of the main examples we have in mind, has a spherical
model category structure which is not strict (\S \ref{dssmc}). This
explains the significance of Remark \ref{ssg} in our context\vsm.

\item The previous example extends to any category $\hC$ of (possibly graded) 
universal algebras with an underlying group structure \ -- \ such as
rings, $R$-modules,  associative algebras, Lie algebras, and so on \ -- \ 
so that $\C$ is corepresented by a $\fG$-theory $\Theta$, in the
language of \cite[\S 4]{BPescF}. Here $\A$ consists
of free monogenic algebras (one for each isomorphism class), and
thus once more \ $\PCAlg\approx\C$. 
\item We can iterate the process by taking $\G$ for $\hC$, and letting \ 
$\hA:=\{\gS{n}\}_{n=1}^{\infty}$ \ (\S \ref{esmod}). We thus obtain
a resolution model category structure on \ $s\G$ \ (or equivalently, on
the category of simplicial spaces). 

In this case the homotopy groups \ $\pi^{s\G}_{k,n}X$, \ denoted
briefly by \ $\pin_{n}X$, \ are the ``bigraded groups'' of
\cite{DKStB}, and Proposition 5.8 there shows that, for a fibrant
simplicial space \ $X\in s\G$, \ we have \ 
$\pia{n}{X}\cong\pi_{0}\,\map(A\otimes S^{n},X)$\vsm.
\item If $\C$ is a resolution model category and $I$ is some small
category, the category \ $\C^{I}$ \ of $I$-diagrams in $\C$ also has
a resolution model category structure, in which the models consists of
all free $I$-diagrams \ $F[A,i]$ \ for $i\in\Obj I$ \ and \ 
$A\in\A$, \ where \ $F[A,i](j):=\coprod_{\Hom_{I}(i,j)}~A$. \ 
See \cite[\S 1]{BJTurR}).
\end{enumerate}

\begin{remark}\label{rbgg}\stepcounter{subsection}
In all these examples, if \ $Y\in \C=s\hC$, \ is fibrant, then for
each \ $n\geq 0$ \ we have an exact sequence: 
%
\setcounter{equation}{\value{thm}}\stepcounter{subsection}
\begin{equation}\label{ethree}
\pis^{\hC}C_{n+1}X~\xra{(\bd^{n+1})_{\#}}~
\pis^{\hC}Z_{n}X~\xra{\hat{\vartheta}_{n}}~ 
\pinC{n}{Y}\to 0.
\end{equation}
\setcounter{thm}{\value{equation}}
\end{remark}

%
%
\sect{Constructions in resolution model categories}
\label{csp}

Not all spherical model categories are resolution model categories 
(see \S \ref{esmc}), but all known examples appear to be Quillen
equivalent to such. Conversely, the examples of resolution model categories \ 
$\lra{\C=s\hC;\A}$ \ we are interested in are spherical (though this
does not hold in general \ -- \ see \S \ref{snsrmc} below). We briefly
indicate why this is so. 

\subsection{Postnikov sections}
\label{sps}\stepcounter{thm}

Given \ $Y\in s\hC$, for each \ $n\geq 0$ \ define \ $Y\q{n}\in s\hC$ \ 
by setting \ $Y\q{n}_{k}:= Y_{k}$ \ for \ $k\leq n+1$ \ and \ 
$Y\q{n}_{k}:= M_{k}(Y\q{n})$ \ (\S \ref{dmat}) for \ $k\geq n+2$. \ 
Note that for any \ $X\in s\C$, \ $M_{k}X$ \ depends only on $X$
through dimension \ $(k-1)$, \ so this definition is valid inductively.
Denote the obvious maps by \ $r\q{n}:Y\to Y\q{n}$ \ and \ 
$p\q{n}:Y\q{n+1}\to Y\q{n}$ \ (see \cite[\S 1.2]{DKanO}).

Now for any \ $X\in s\hC$, \ choose a functorial fibrant replacement
$Y$, and set \ $P_{n}X := Y\q{n}$, \ with \ 
$\varphi\q{n}:X\to P_{n}X$ \ defined to be the composite of \ $r\q{n}$ \ 
with the trivial cofibration \ $i:X\to Y$, \ and \ 
$p\q{n}:P_{n+1}X\to P_{n}X$ \ defined as above. 

\begin{remark}\label{rcsk}\stepcounter{subsection}
The functor \ $-\q{n}:\C\to\C$ \ is right adjoint to the \
$(n+1)$-skeleton functor \ $\sk{n+1}$, \ so \ $P_{n}X$ \ depends only
on \ $\sk{n+1}X$, \ even if $X$ is not fibrant. If $X$ is fibrant, we can find \ 
$Y\simeq P_{n}X$ \ with \ $\sk{n+1}Y=\sk{n+1}X$.
\end{remark}

%
%
\begin{fact}\label{fone}\stepcounter{subsection}
In each of the examples of \S \ref{ermc}(a-d), the tower:
$$
X\to\dotsc \to P_{n+1}X\xra{p\q{n}} P_{n}X\to\dotsc \to P_{0}X
$$
is a functorial Postnikov tower for \ $\C=s\hC$ \ with respect to $\A$ \ 
(\S \ref{dfpt}).  
\end{fact}

\begin{proof}
From \S \ref{sreed} and \S \ref{sps} it follows that if  \ $Y\in
s\hC$ \ is fibrant, then \ so is each \ $Y\q{n}$, \ and  
for each $n$, \ $Y\q{n+1}\to Y\q{n}$ \ is a fibration, \ 
$Z_{k}Y\q{n}=0$ \ and \ $C_{k}Y\q{n}\xra{\bd}Z_{k-1}Y$ \ is an
isomorphism for \ $k\geq n+2$. \ The claim then follows from \eqref{ethree}.
\end{proof}

%
%
\begin{fact}\label{ftwo}\stepcounter{subsection}
In each of the examples of \S \ref{ermc}(a-d), \ there is a
classifying object \ $\BL$ \ for any \PAa\ \ $\Lambda$, \ and it is
unique up to homotopy. 
\end{fact}

\begin{proof}
In the algebraic cases of \S \ref{ermc}(a-b), \ we may take \ $\BL$ \
to be (a cofibrant model for) the constant simplicial object on
$\Lambda$. For simplicial spaces, \ $\BL$ \ may be constructed as for
topological spaces, using generators and relations (see \cite[\S 8.9]{BDGoeR}). 
The extension to the diagram case of \S \ref{ermc}(d) is objectwise.
\end{proof}

%
%
\begin{fact}\label{fthree}\stepcounter{subsection}
In each of the examples of \S \ref{ermc}(a-d), \ for each \ $n\geq 1$ \ 
there is an $n$-dimensional $M$-Eilenberg-Mac~Lane object \ 
$\EK{}{}{M}{n}$ \ for any abelian \PAa\ $M$, \ and there is an
$n$-dimensional extended $M$-Eilenberg-Mac~Lane object \ 
$\EL{M}{n}$ \ for any \PAa\ $\Lambda$ and module $M$ over
$\Lambda$. \ Each of these is unique up to homotopy. 
\end{fact}

\begin{proof}
In the algebraic cases of \S \ref{ermc}(a-b), \ we may take \
$\EK{}{}{M}{n}$ \ to be the iterated Eilenberg-Mac Lane construction
$\bar{W}$ on \ $BM$, \ while \ $\EL{M}{n}$ \ is a semi-direct product \ 
$\EK{}{}{M}{n}\ltimes\BL$ (see \cite[Prop.\ 2.2]{BDGoeR}).
For simplicial spaces, use the explicit construction of \cite[\S 8.9]{BDGoeR} \
The extension to the diagram case is again objectwise.
\end{proof}

%
%
\begin{fact}\label{ffour}\stepcounter{subsection}
In each of the examples of \S \ref{ermc}(a-d), \ for each \ $n\geq 1$ \ 
and \ $X\in\C$, \ $\pinC{n}X$ \ has a natural structure of a 
module over \ $\pinC{0}X$.
\end{fact}

\begin{proof}
Note that by \cite[II,1,(6)]{QuiH} we have \ 
$\map(A\otimes S^{n},X)\cong\map_{\Ss}(S^{n},\map(A,X))$ \ 
(unpointed maps), \ so \ $\pinC{n}X\to\pinC{0}X$ \ associates to each \ 
$f:A\otimes S^{n}\to X$ \ its component in \ $\map(A,X)$. \ This defines 
an \emph{abelian} algebra over \ $\pinC{0}X$ \ by \cite[Prop.\ 6.26]{BPescF}).
\end{proof}

%
%
\begin{fact}\label{ffive}\stepcounter{subsection}
In each of the examples of \S \ref{ermc}(a-d), \ for each \ $X\in s\C$, \ 
$\Lambda:=\pinC{0}X$ \ and \ $n\geq 1$, \ the commutative square
obtained by applying the functor \ $P_{n+2}$ \ to the pushout diagram:
$$
\xymatrix@C=40pt{\ar @{} [dr] |>>>>>{\framebox{\scriptsize{PO}}}
P_{n+1}X  \ar[r]^-{p\q{n+1}} 
\ar[d] & P_{n}X \ar[d]^{k_{n}}\\ \BL \ar[r] & Y
}
$$
is an $n$-th $k$-invariant square (Def.\ \ref{dki}) \ -- \ that is, \ 
$P_{n+2}Y\simeq\EL{\pia{n+1}X}{n+2}$.
\end{fact}

\begin{proof}
See \cite[\S 5]{BDGoeR}.
\end{proof}

We may summarize these facts in the following:

%
%
\begin{thm}\label{tone}\stepcounter{subsection}
The following resolution model categories (cf.\ \S \ref{ermc})
are strict spherical model categories:
\begin{enumerate}
\renewcommand{\labelenumi}{\roman{enumi}.\ }
\item The category \ $\C=s\hy{\Theta}{\Seta}$ \ of simplicial
  $\Theta$-algebras for any $\fG$-theory $\Theta$,
  with $\hA$ consisting of monogenic free $\Theta$-algebras;
\item In particular, the category \ $\C=\G$ \ of simplicial groups, with \
  $\hA=\{\bZ\}$; 
\item The category \ $s\G$ \ of bisimplicial groups (``simplicial
  spaces''), with \ $\hA=\{\gS{1}\otimes S^{k}\}_{k=0}^{\infty}$.
\item The category \ $\C^{I}$ \ of $I$-diagrams in a strict spherical
  model category $\C$.
\end{enumerate}
\end{thm}

%
%
\begin{thm}\label{ttwo}\stepcounter{subsection}
The following are spherical model categories (which are not strict):
\begin{enumerate}
\renewcommand{\labelenumi}{\roman{enumi}.\ }
\item The category \ $\Sa$ \ of pointed simplicial sets, with \ $\A=\{S^{1}\}$; 
\item The category \ $\Ta$ \ of pointed topological spaces, with \ 
      $\A=\{\bS{1}\}$; 
\item The category \ $s\Ta$ \ of simplicial pointed topological spaces, with \ 
      $\hA=\{\bS{1}\otimes S^{k}\}_{k=1}^{\infty}$.
\end{enumerate}
\end{thm}

\subsection{Non-spherical model categories}
\label{snsrmc}\stepcounter{thm}

Consider the trivial model category structure on \ $\hC=\Gp$, \ with \ 
$\hA:=\{A=\bZ/p\}$ \ (for $p$ an odd prime). \ This defines a
resolution model category structure on $\G$ \ -- \ or equivalently, on \ 
$\Ta$ \ (see Remark \ref{rmods}). \ Note that \ $-\otimes S^{n}$ \ 
corresponds to suspension of simplicial sets, not
simplicial abelian group, so the model \ $A\otimes S^{n}\in\G$ \
corresponds to the $n$-dimensional mod $p$ Moore space \ 
$\bS{n-1}\cup_{p}\be{n}$. 

Thus \ $\pia{k}X:=[A\otimes S^{k},X]$ \ is by definition the $k$-th
\emph{mod $p$ homotopy group} of $X$ \ -- \  denoted by \ $\pi_{k}(X;\bZ/p)$ \ 
in \cite[Def.\ 1.2]{NeiP} \ -- \  which fits into a short exact sequence:
%
\setcounter{equation}{\value{thm}}\stepcounter{subsection}
\begin{equation}\label{efour}
0\to\pi_{k}X\otimes\bZ/p \to
        \pi_{k}(X;\bZ/p)\to\Tor^{\bZ}_{1}(\pi_{k-1}X,\bZ/p)\to 0
\end{equation}
\setcounter{thm}{\value{equation}}
\noindent for \ $k\geq 2$ \ (see \cite[Prop.\ 1.4]{NeiP}). \ In particular, 
for \ $Y:=A\otimes S^{n}$ \ ($n\geq 4$) \ we have
$$
\pi_{i}(Y;\bZ/p)=\begin{cases}\bZ/p & \text{for \ }i=n-1,n,\\
             0     & \text{for \ } 2\leq i<n-1 \text{\ or \ } i=n+1,
\end{cases}
$$
\noindent with the two non-trivial groups connected by a Bockstein 
(cf.\ \cite[\S 1]{NeiP})\vsm.

However, the resolution model category structure on $\G$ determined by $\A$ is 
not spherical: \ if it were, in particular there would be Postnikov functors \
$P_{k}=P_{k}^{\A}$ \ for all \ $k\geq 1$ \ (Def.\ \ref{dfpt}). \ 
From \eqref{efour} we see that, disregarding torsion prime to $p$,  
because of the Bockstein we must have \ 
$P_{n-1}Y\simeq E(\bZ,n-1)$ \ and \ $P_{n}Y\simeq E(\bZ/p,n-1)$ \ 
(for \ $Y=S^{n-1}\cup_{p} e^{n}$). \ But then there is no
non-trivial map \ $P_{n}Y\to P_{n-1}Y$. 

\subsection{Cohomology in spherical model categories}
\label{scpa}\stepcounter{thm}

Note that the  $k$-invariants of a simplicial object
actually take value in cohomology groups, as expected:

%
%
\begin{prop}\label{ptwo}\stepcounter{subsection}
For each \PAa\ $\Lambda$ and module $M$ over $\Lambda$, the functors \
$D^{n}:\C/\BL\to\Abgp$ \ \textup{($n>0$)}, \ defined \ 
$D^{n}(X):=[X,\EL{M}{n}]_{\BL}$, \ 
are \emph{cohomology functors} on $\C$ \ -- \ that is, they 
are homotopy invariant, take arbitrary coproducts to products, vanish 
on the spherical models \ $\Sigma^{n}A$, \ except in degree $n$, 
and have Mayer-Vietoris sequences for homotopy pushouts.
\end{prop}

We therefore denote \ $[X,\EL{M}{n}]_{\BL}$ \ by \ $\HL{n}{X}{M}$.

\begin{proof}
See \cite[Thm.\ 7.14]{BPescF}.
\end{proof}

Fact \ref{fthree} then follows from Brown Representability, since \
$\EL{M}{n}$ \ represents the $n$-th Andr\'{e}-Quillen cohomology
group in $\C$; see  \cite[\S 6.7]{BDGoeR} and \cite[\S 4]{BlaG}.  

%
%
\sect{Spherical functors}
\label{csf}

Our objective is to study functors between model categories, and
investigate the extent to which they induce an equivalence of homotopy
categories. Our methods work only for functors between spherical model
categories which take models to models, in the following sense:

\begin{defn}\label{dsfunc}\stepcounter{subsection}
Let \ $\lra{\C;\A}$ \ and \ $\lra{\D;\B}$ \ be two spherical model
categories. A functor \ $T:\C\to\D$ \ is called \emph{spherical} if 
\begin{enumerate}
\renewcommand{\labelenumi}{\roman{enumi}.\ }
\item $T$ defines a bijection \ $\A\to\B$;
\item $T\rest{\PiA}$ \ preserves coproducts and suspensions;
\item $T$ induces an equivalence of categories \ $\PCAlg\approx\PDAlg$ \ 
(in fact, it suffices that \ $\PDAlg$ \ be a full subcategory of \ $\PCAlg$).
\end{enumerate}
\end{defn}

\subsection{Examples of spherical functors}
\label{ssfunc}\stepcounter{thm}

In the cases we shall be considering (those mentioned in the
introduction), $\C$ and $\D$ will be strict spherical resolution model
categories, with \ $\C=s\hC$ \ and \ $\D=s\hD$, \ and $T$ will 
be prolonged from a functor \ $\hT:\hC\to\hD$\vsm . 

The four examples:

\begin{enumerate}
\renewcommand{\labelenumi}{(\alph{enumi})\ }
\item For \ $\lra{\hC;\hA}=\lra{\Gp;\{\bZ\}}$ \ and \  
$\lra{\hD,\hB}=\lra{\Abgp;\{\bZ\}}$, \ let \ $\hT=\Ab:\Gp\to\Abgp$ \
be the abelianization functor.

Here \ $\C=s\hC=\G$, \ so \ $\ho\C$ \ is equivalent to the homotopy
category of pointed connected topological spaces (\S \ref{ssg}), 
while \ $\D=s\hD$, \ the category of simplicial abelian groups,
is equivalent to the category of chain complexes under the
Dold-Kan correspondence (see \cite[\S 1]{DolH}). Thus \ $T:\C\to\D$ \
represents the singular chain complex functor \ $C_{\ast}:\Ta\to\Ch$.

Note that \ $\PCAlg=\Gp$, \ while \ $\PDAlg=\Abgp$, \ in this case, so 
strictly speaking $T$ does not induce an equivalence of categories. \
But since \ $\Abgp$ \ is a full subcategory of \ $\Gp$, \ we can in
fact think of \ $\pin$ \ as taking values in groups\vsm .
\item For \ $\lra{\hC;\hA}=\lra{\Gp;\{\bZ\}}$ \ and \  
$\lra{\hD,\hB}=\lra{\Hopf;\{H\}}$, \ where \ $\Hopf$ \ is the
category of complete Hopf algebras over $\bQ$, \ $H$ is the
monogenic free object in this category, let \ $\hQ:\Gp\to\Hopf$ \ 
be the functor which associates to a group $G$ the completion of the
group ring \ $\bQ[G]$ \ by powers of the augmentation ideal. 

Again, \ $\C=s\hC$ \ is a model category for connected 
topological spaces, while \ $\D=s\hD$ \ is a model category for the
rational simply-connected spaces (see \cite{QuiR}); \ $Q$ (when
restricted to connected simplicial groups) represents the
rationalization functor. \ Once more, \ $\PCAlg=\Gp$, \ while \
$\PDAlg$ \ is the subcategory of vector spaces over $\bQ$\vsm. 
\item For \ $\lra{\hC,\hA}=\lra{\Seta;\{S^{0}\}}$ \ 
(so that \ $\lra{\C,\A}=\lra{\Ss;\{S^{1}\}}$, \ by Remark
\ref{rmods}), and \ $\lra{\hC;\hA}=\lra{\Gp;\{\bZ\}}$, \ let \
$\hF:\Seta\to\Gp$ \ be the free group functor.

Again, we think of both \ $\C=s\hC=\G$ \ and \ $\D=s\hD=\Sa$ \ as
model categories for pointed topological spaces, (under the respective
equivalences of \S \ref{ssg}) \ -- \ so $F$ here represents 
the suspension functor \ $\Sigma:\Ta\to\Ta$ \ (rather than \
$\Omega\Sigma$, \ as one might think at first glance)\vsm.
\item For \ $\lra{\hC;\hA}=\lra{\G;\{\gS{k}\}_{k=0}^{\infty}}$ \ and \  
$\lra{\hD,\hB}=\lra{\PAlg;\{\pis\gS{k}\}_{k=0}^{\infty}}$, \ let \ 
$\hP:\G\to\PAlg$ \ be the graded homotopy group functor \ $X\mapsto\pis X$. \ 
Here \ $\C=s\G$ \ is a model category for simplicial spaces\vsm.
\end{enumerate}

%
%
\begin{thm}\label{tthree}\stepcounter{subsection}
Let \ $\lra{\C;\A}$ \ and \ $\lra{\D;\B}$ \ be spherical model
categories,  and let \ $T:\C\to\D$ \ be a spherical functor. Then for
each \ $X\in\C$ \ and \ $A\in\A$ \ there is a natural long exact
sequence of \PCa s: 
%
\setcounter{equation}{\value{thm}}\stepcounter{subsection}
\begin{equation}\label{efive}
\dotsc\to \Gamma^{X}_{\alpha,n}X~\xra{s_{n}^{X}}~\pia{n}
X~\xra{h_{n}^{X}}~\pinC{T_{\ast}(\alpha),n}
TX~\xra{\partial_{n}^{X}}~\Gamma^{T}_{\alpha,n-1}X\dotsc~. 
\end{equation}
\setcounter{thm}{\value{equation}}
\end{thm}

We call \eqref{efive} the \emph{comparison exact sequence} for
$T$. Compare \cite[V, (5.4)]{BauCF}.

\begin{proof}
If \ $\tX\to X$ \ is a functorial fibrant replacement, the functor 
$T$ induces a natural transformation \ 
$\tau:\map_{\C}(A,\tX)\to \map_{\D}(TA,\widehat{T\tX})$, \ which
we may functorially change to a fibration of simplicial sets, 
with fiber \ $F(X)$. \ Setting \ $\Gamma^{T}_{\alpha,n}:=\pi_{n}F(X)$, \ 
the corresponding long exact sequence in homotopy is
\eqref{efive}. 

Note that the map \ $h^{X}_{n}=h^{X}$ \ is also natural in the
variable \ $A$, \ so the graded map \ $h^{X}_{\ast}:\pinC{n}X\to\pinD{n}TX$ \ 
is a morphism of \PCa s (i.e., \ $\PihA$-algebras).
\end{proof}

\subsection{Applications of Theorem \protect{\ref{tthree}}}
\label{satt}\stepcounter{thm}

The Theorem is not very useful in this generality. However, in all 
the examples of \S \ref{ssfunc}, we obtain
interesting (though mostly known) exact sequences\vsm :

\begin{enumerate}
\renewcommand{\labelenumi}{(\alph{enumi})\ }
\item For \ $\hT=\Ab:\Gp\to\Abgp$ \ the abelianization functor, where \ 
$T:\G\to s\Abgp$ \ represents the singular chain complex functor \ 
$C_{\ast}:\Ta\to\Ch$ \ (cf.\ \S \ref{ssfunc}(a)), the sequence 
\eqref{efive} is the ``certain exact sequence'' of J.H.C. Whitehead:
%
\setcounter{equation}{\value{thm}}\stepcounter{subsection}
\begin{equation}\label{etwentyfive}
\dotsc\to \Gamma_{n} X \to\pi_{n} X~\xra{h_{n}}~ H_{n}(X;\bZ)\to
\Gamma_{n-1} X \dotsc
\end{equation}
\setcounter{thm}{\value{equation}}

(See \cite{JWhC}).  In particular, the third term in this sequence, \
$\Gamma^{\A}_{n}(X)$, \ is simply the $n$-th homotopy group of the
commutator subgroup of \ $GX$\vsm. 
\item For \ $Q:\G\to s\Hopf$ \ of \S \ref{ssfunc}(b), representing 
the rationalization functor, \ we obtain a long exact sequence
relating the integral and rational homotopy groups of a
simply-connected space $X$. The third term in \eqref{efive} may be
described in terms of the torsion subgroup of \ $\pis X$ \ together
with \ $\pis X\otimes\bQ/\bZ$\vsm. 

\item The free group functor \ $\hF:\Seta\to\Gp$ \ of \S
\ref{ssfunc}(c) represents the suspension \ $\Sigma:\Ta\to\Ta$, \ 
and indeed for \ $K\in\Sa$ \ the map \ $h^{K}$, \ which is the composite:

\begin{equation*}
\begin{split}
\pi_{n}K~=~\pi_{0}\map_{\Sa}(S^{n},K)~&~
\longrightarrow~\pi_{0}\,\map_{\G}(FS^{n},FK)\\
\xra{\cong}~& \pi_{0}\,\map_{\Sa}(\Sigma S^{n},\Sigma K)~=~\pi_{n+1}\Sigma K~,
\end{split}
\end{equation*}
\noindent is the suspension homomorphism, so \eqref{efive} is a
generalized EHP sequence (cf.\ \cite{BauR,GanG,NomE})\vsm. 
\item For \ $\pis:s\G\to s\PAlg$ \ as in \S \ref{ssfunc}(d), it
turns out that for any simplicial space \ $X\in s\G$, \ the induced map \
$h^{X}_{n}$ \ is the ``Hurewicz homomorphism'' \ 
$h_{n}:\pin_{n}X\to\pi_{n}\pis X$ \ of \cite[7.1]{DKStB}, \ 
while \ $\Gamma^{T}_{n}X$ \ is just \ $\Omega\pin_{n-1}X$ \ -- \ that is, \ 
$\Gamma^{T}_{i,n}X=\pin_{i+1,n-1}X$ \ for each $i$. Thus \ 
\eqref{efive} \ is the \emph{spiral long exact sequence}:
%
\setcounter{equation}{\value{thm}}\stepcounter{subsection}
\begin{equation}\label{efifteen}
\dotsc \pi_{n+1}\pis X \xra{\partial^{\star}_{n+1}} \Omega\pin_{n-1}X
\xra{s_{n}} \pin_{n}X \xra{h_{n}} \pi_{n}\pis X\to
\dotsb \pin_{0}X\xra{h_{0}}\pi_{0}\pis X \to 0
\end{equation}
\setcounter{thm}{\value{equation}}
\noindent  of \cite[8.1]{DKStB}. Of course, \ $\pin_{-1}{X}=0$, \ 
so \ $h_{0}$ \ is an isomorphism\vsm .
\end{enumerate}

Note that for \ $T:\C\to\D$ as above, the homotopy groups \ $\pinD{n}TX$ \ for
any \ $X\in\C=s\hC$ \ may be computed using the Moore chains \
$C_{\ast}TX$ \ as in \S \ref{dnc}; each  \ $\pinD{n}TX$ \ is a \ $\PiD$-algebra,
abelian for \ $n\geq 1$.  

\subsection{Explicit construction of the spiral exact sequence}
\label{sses}\stepcounter{thm}

It may be helpful to inspect in detail the construction of last long
exact sequence, since it is perhaps the least familar of the four.
Specificializing to \ $\hC=\G$ \ and \ $T=\pis$, \ we have:

%
%
\begin{lemma}\label{lone}\stepcounter{subsection}
For fibrant \ $X\in\C$, \ the inclusion \ $\iota:C_{n}X\hra X_{n}$ \ 
induces an isomorphism \ $\iota_{\star}:\pis C_{n}X\cong C_{n}(\pis X)$ \ 
for each \ $n\geq 0$. 
\end{lemma}

\begin{proof}
See \cite[Prop.\ 2.11]{BlaCW}.
\end{proof}

Together with \eqref{ethree}, this yields a commuting diagram:
%
\setcounter{equation}{\value{thm}}\stepcounter{subsection}
\begin{equation}\label{etwentyseven}
\xymatrix{
\pis C_{n+1}X \rto^{(\bd)_{\#}} \dto_{\iota_{\star}}^{\cong} & 
\pis Z_{n}X \ar@{->>}[r]^{\hat{\vartheta}_{n}} \dto^{\hat{\iota}_{\star}} & 
\pin_{n}X \ar@{.>}[d]^{h_{n}} \\
C_{n+1}(\pis X) \rto^{d_{0}^{\pis X}} & 
Z_{n}(\pis X) \ar@{->>}[r]^{\vartheta_{n}} & \pin_{n}\pis X
}
\end{equation}
\setcounter{thm}{\value{equation}}
\noindent which defines the dotted morphism of \Pa s \ 
$h_{n}:\pin_{n}X\to\pi_{n}(\pis X)$. \ Note that for \ $n=0$ \ the map \ 
$\hat{\iota}_{\star}$ \ is an isomorphism, so $h$ is, too.

If  \ $X\in s\G$ \ is fibrant, applying \ $\pis$ \ to the fibration sequence 
\eqref{etwo} yields a long exact sequence, with connecting homomorphism \ 
$\partial_{n}:\Omega\pis Z_{n} X=\pis\Omega Z_{n} X \to\pis Z_{n+1} X $; \ 
\eqref{ethree} then implies that
%
\setcounter{equation}{\value{thm}}\stepcounter{subsection}
\begin{equation}\label{efourteen}
\Omega\pin_{n}X=\Omega\Cok(\bd^{n+1})_{\#}\cong\Image\partial_{n}
\cong \Ker(j^{X}_{n+1})_{\#}\subseteq\pis Z_{n+1}X,
\end{equation}
\setcounter{thm}{\value{equation}}
\noindent and the map \ $s_{n+1}:\Omega\pin_{n}X\to\pin_{n+1}X$ \ in
\eqref{efour} is then obtained by composing the inclusion \ 
$\Ker(j^{X}_{n+1})_{\#}\hra\pis Z_{n+1}X$ \ with the quotient map \ 
$\hat{\vartheta}_{n+1}:\pis Z_{n+1}X\to\pin_{n+1}X$ \ of \eqref{ethree}.

Similarly, \ $h_{n}:\pin_{n}X\to\pin_{n}\pis X$ \ is induced by the 
inclusion \ 
$(j^{X}_{n})_{\#}:\pis Z_{n}X\to Z_{n}\pis X\subseteq C_{n}\pis X$, \ and \ 
$\partial^{\star}_{n+2}:\pin_{n+2}\pis X\to\Omega\pin_{n}X$ \ is induced by \ 
the composite \ 
$$
Z_{n+2}\pis X\subseteq C_{n+2}\pis X\cong \pis C_{n+2}X
\xra{(\bd^{n+2})_{\#}}Z_{n+1}\pis X,
$$
\noindent which actually lands in \ 
$\Ker(j^{X}_{n+1})_{\#}\cong \Omega\pin_{n}X$ \ by the exactness of the 
long exact sequence for the fibration.

Moreover, for each \ $n\geq 0$, \ \eqref{etwentyseven} may be extended
(after rotating by $90^{\circ}$) \ to a commutative diagram with exact
rows and columns: 
$$
\xymatrix{
& 0 \dto & 0\dto & 0\dto & & \\ 
0 \rto & \Ker s_{n} \ar@{|->}[r] \dto & B_{n+1}X \dto
\ar@{->>}[r]^>>>>>{(j_{n})_{\ast}} & B_{n+1}\pis X_{n+2}\rto\dto & 0\dto & \\
0 \rto & \Omega\pin_{n-1}X \ar@{|->}[r] \dto & 
\pis Z_{n}X \dto^{\hat{\vartheta}_{n}} \rto^{(j_{n})_{\ast}} & 
Z_{n}\pis X \ar@{->>}[r] \dto^{\vartheta_{n}} & \Cok h_{n} \rto \dto^{=} & 0 \\
0 \rto & \Ker h_{n} \ar@{|->}[r] \dto & \pin_{n}X \dto \rto^{h_{n}} & 
\pi_{n}\pis X \ar@{->>}[r] \dto & \Cok h_{n}\rto\dto & 0 \\
 & 0 & 0 & 0 & 0 &
}
$$
\noindent in which \ 
$B_{n+1} X:=\Image(\bd^{X_{n+2}})_{\#}\subseteq\pis Z_{n} X$ \ and \  
$B_{n+1}\pis X_{n+2}:=\Image\bd^{\pis X_{n+2}}$ \ are the respective
boundary objects.

The maps \ $\partial^{\star}_{n+1}$, \ $s_{n}$, \ and \ $h_{n}$, \ 
as defined above, form the spiral long exact sequence.

\subsection{Inverse spherical functors}
\label{sisf}\stepcounter{thm}

We may sometimes be interested in functors between spherical model
categories which are not quite spherical. Thus, if \ 
$T:\lra{\C;\A}\to\lra{\D;\B}$ \ is a spherical functor as in \S
\ref{dsfunc}, a functor \ $V:\D\to\C$ \ equipped with a natural
transformation \ $\vartheta:\Id_{\C}\to VT$ \ is called an
\emph{inverse spherical functor} to $T$. 

\begin{example}\label{eisf}\stepcounter{subsection}
For the free group functor \ $F:\Seta\to\Gp$ \ of \S \ref{ssfunc}(c),
the forgetful functor \ $\hU:\Gp\to\Seta$ \ (right adjoint to $F$)
with the adjunction counit \ $\eta:\Id\to UF$ \ as the natural
transformation $\vartheta$, yields the inverse spherical functor \ 
$U:\G\to\Sa$. \ Here we do not think of $\G$ as a model for $\Ta$ \ -- \ 
rather, $U$ represents the forgetful functor from loop spaces
(topological groups) to spaces.

Similarly, the adjoint to the abelianization functor \ $\Ab:\Gp\to\Abgp$ \ 
is the inclusion \ $\hat{I}:\Abgp\to\Gp$, \ and the corresponding functor \ 
$I:s\Abgp\to\G$ \ represents the factorization of the Dold-Thom infinite
symmetric product functor \ $SP^{\infty}:\Ta\to\Ta$ \ 
through \ $\Ch$.
\end{example}

%
%
\begin{prop}\label{pfour}\stepcounter{subsection}
If \ $V:\D\to\C$ \ is an inverse spherical functor to $T$,
then for each \ $Y\in\D$ \ and \ $B\in\B$ \ there is a natural 
long exact sequence:
%
\setcounter{equation}{\value{thm}}\stepcounter{subsection}
\begin{equation}\label{esix}
\dotsc\to \Delta^{V}_{B,n}Y\to\pinD{B,n} Y\xra{V_{\#}} 
\pinC{V_{\ast}(B),n} VY\to \Delta^{V}_{B,n-1}Y\dotsc
\end{equation}
\setcounter{thm}{\value{equation}}
\end{prop}

\begin{proof}
If $V$ is an inverse spherical functor, because \ $T\rest{\A}$ \ is a
bijection onto $\B$, there is an \ $A\in\A$ \ such that \ 
$B=TA$. \ As before, $V$ induces a natural transformation \ 
$\nu:\map_{\D}(B,\tY)\to \map_{\D}(VB,\widehat{V\tY})$ \ and the
natural transformation \ $\vartheta:A\to VTA$ \ yields \ 
$\vartheta^{\#}:\map_{\D}(VTA,\widehat{V\tY})\to
\map_{\D}(A,\widehat{V\tY})$ \ so we get a composite map \ 
$\map_{\D}(B,\tY)\to\map_{\D}(A,\widehat{V\tY})$, \ with homotopy
fiber \ $E(Y)$. \ If we let \ $\Delta^{V}_{\beta,n}Y:=\pi_{n} E(Y)$, \ 
the fibration long exact sequence is \eqref{esix}. 
\end{proof}

\begin{remark}\label{rlesisf}\stepcounter{subsection}
Note that in contradistinction to Theorem \ref{tthree}, \ 
$V_{\#}$ \ of \eqref{esix} need not respect any operations, since we
only have a bijection \ $T\rest{\A}:\A\to\B$, \ not a functor. 

For \ $U:\G\to\Sa$ \ as in \S \ref{eisf}, we may assume \ $X\in\G$ \
is of the form \ $X\simeq GK$ \ for \ $K\in\Sa$, \ and then \ 
$V_{\#}$ \ is the identity:
\setcounter{equation}{\value{thm}}\stepcounter{subsection}
\begin{equation}\label{etwentythree}
\begin{split}
\pi_{n}K~=~\pin_{n}X&=~\pi_{0}\,\map_{\G}(FS^{n-1},GK)~\to~
\pi_{0}\,\map_{\Sa}(UFS^{n-1},UGK)\\
& \xra{\eta^{\#}}~\pi_{0}\,\map_{\Sa}(S^{n-1},UGK)~=~\pi_{n}K~,
\end{split}
\end{equation}
\setcounter{thm}{\value{equation}}
\noindent so \eqref{esix} is not interesting in this case. 
\end{remark}

%
%
\sect{Comparing Postnikov systems}
\label{ccps}

The basic problem under consideration in this paper may be formulated as
follows\vsm : 

\subsection*{Question}
Given a spherical functor \ $T:\lra{\C;\A}\to\lra{\D;\B}$ \
and an object \ $G\in\D$, \ what are the different objects \ $X\in\C$ \ 
(up to homotopy) such that \ $TX\simeq G$\vsm ?

As shown in the previous section, such a pair \ 
$\lra{X,G}$ \ must be connected by a comparison exact sequence. Thus,
in order to reconstruct $X$ from $G$, we first try to determine \ 
$\pinC{\ast}X$, \ and its relation to \ $\pinD{\ast}G$. 

In order to proceed further, we must make an additional assumption on
$T$, contained in the following:

\begin{defn}\label{dspecial}\stepcounter{subsection}
A spherical (or inverse spherical) functor \ $T:\C\to\D$ \ is called
\emph{special} if:

\begin{enumerate}
\renewcommand{\labelenumi}{\roman{enumi}.\ }
\item $\C=s\hC$ \ and \ $\D=s\hD$ \ are spherical resolution model
  categories, and $T$ is prolonged from a functor \ $\hT:\hC\to\hD$. \ 
\item For any \PAa\ $\Lambda$ and module $M$ over $\Lambda$, \ $T$ induces a 
homomorphism of (graded) groups \ $\phi_{T}:\Lambda\to\pinD{0}T\BCL$. \ 
\item This \ $\phi_{T}$ \ induces a functor \ 
$\hat{T}:\RM{\Lambda}\to\RM{\phi_{T}\Lambda}$ \ which is an
isomorphism on $\Lambda$-modules (see Remark \ref{rmodule}).
\item For each \ $n\geq 1$ \ and $n$-dimensional extended
$M$-Eilenberg-Mac~Lane object \ $E=\ECL{M}{n}$, \ there is a natural
isomorphism \ $\pinD{n}TE\cong M$ \ which respects $\hat{T}$ in the
obvious sense.
\item The natural map
%
\setcounter{equation}{\value{thm}}\stepcounter{subsection}
\begin{equation}\label{etwentyfour}
[X,\ECL{M}{n}]_{\BCL}\to [TX,E\sp{\hat{T}\Lambda}\sb{\D}(M,n)]_{B_{\D}\hat{T}L}~,
\end{equation}
\setcounter{thm}{\value{equation}}
\noindent defined by composition with the projection
$$
\rho:T\ECL{M}{n}\to P_{n}T\ECL{M}{n}~=~E\sp{\hat{T}\Lambda}\sb{\D}(M,n)~, 
$$ 
is an isomorphism.
\end{enumerate}
\end{defn}

\begin{example}\label{especial}\stepcounter{subsection}
All the functors we have considered hitherto, except for the rationalization
functor \ $Q:\G\to s\Hopf$ \ of \S \ref{ssfunc}(b), are special\vsm:

\begin{enumerate}
\renewcommand{\labelenumi}{(\alph{enumi})\ }
\item For the singular chain functor \ $T:\G\to s\Abgp$, induced by 
abelianization, this follows from the Hurewicz Theorem (recall that \
$\pinC{0}X$ \ is the fundamental group, in our indexing for \ $X\in\G$).
\item For the suspension \ $\Sigma:\Ta\to\Ta$, \ induced by 
the free group functor \ $F:\Seta\to\Gp$, \ this follows (in the
simply connected case) from the Freudenthal Suspension Theorem.
\item For the homotopy groups functor \ $\pis:s\G\to s\PAlg$, \ (i)-(iii)
follow by inspecting the spiral long exact sequence \eqref{efifteen},
while (iv) is \cite[Prop.\ 8.7]{BDGoeR}.
\item For the inverse spherical functor \ $U:\G\to\Sa$ \ of \S \ref{eisf},
induced by the forgetful functor \ $\hU:\Gp\to\Seta$, \ this is
immediate from \eqref{etwentythree}\vsm.
\end{enumerate}
\end{example}

%
%
\begin{lemma}\label{ltwo}\stepcounter{subsection}
Any special spherical functor \ $T:\C\to\D$ \ as above \emph{respects
Postnikov systems} \ -- \ that is, for any \ $X\in\C$ \ 
and \ $n\geq 0$ \ we have:
%
\setcounter{equation}{\value{thm}}\stepcounter{subsection}
\begin{equation}\label{eeight}
P^{\D}_{n}T P^{\C}_{n} X\cong P^{\D}_{n}T X \ - \ 
\end{equation}
\setcounter{thm}{\value{equation}}
so that \ $\pinC{k}T X\cong \pinD{k}TP_{n}X$ \ and \ 
$\Gamma_{k}X\cong \Gamma_{k}P_{n}X$ \ for \ $k\leq n$. 
\end{lemma}

\begin{proof}
This follows from the constructions in \S \ref{sps} and the proof of
Theorem \ref{tthree}.
\end{proof}

\subsection{Postnikov systems and spherical functors}
\label{spssf}\stepcounter{thm}

From now on, assume \ $T:\C\to\D$ \ is a special spherical functor.
Ultimately, for each object \ $G\in\D$, \ we would like find any and all \
$X\in\C$ \ such that \ $TX\simeq G$. \ First, however, we try to
discover what can be said about \ $TX$ \ and its Postnikov systems
for a given \ $X\in\C$. \ Using the comparison exact sequence for $T$
and Lemma \ref{ltwo}, we see that:  
%
\setcounter{equation}{\value{thm}}\stepcounter{subsection}
\begin{equation}\label{enine}
\pinD{k}TP_{n}X\cong\begin{cases}
        \pinD{k}TX & \text{for \ } k\leq n,\\
        \Cok\{h^{X}_{n+1}:\pinC{n+1}X\to\pinD{n+1}T X\} & \text{for \ }
          k=n+1,\\
        \Gamma_{k-1}P_{n}X & \text{for \ }k\geq n+2~.
\end{cases}
\end{equation}
\setcounter{thm}{\value{equation}}

%
%
\begin{fact}\label{fsix}\stepcounter{subsection}
If \ $T:\C\to\D$ \ is a special spherical functor, applying \
$\pinC{n+2}$ \ to the $n$-th $k$-invariant \
$k_{n}:P_{n}X\to \ECL{\pinC{n+1}X}{n+2}$ \  
yields the homomorphism \ $s^{X}_{n+1}:\Gamma_{n+1}X\to\pinC{n+1}X$.
\end{fact}

\begin{proof}
Since $T$ is special, \ $\pinD{n+2}T\ECL{\pin_{n+1}X}{n+2}\cong\pinC{n+1}X$, \ 
and \ $\pinD{n+2}TP_{n}X\cong\Gamma_{n+1}X$ \ from 
\eqref{enine}, so this follows from the naturality of the comparison
exact sequence, applied to the maps in \eqref{ezero}.
\end{proof}

%
%
\begin{lemma}\label{lthree}\stepcounter{subsection}
If \ $T:\C\to\D$ \ is a special spherical functor, for any \ $X\in\C$, 
$$
\xymatrix{
P_{n+1}TP_{n}X \rto \dto & P_{n+1}TP_{n-1}X \dto \\ 
P_{n+1}T\BCL \rto^-{\protect{Tk_{n}}} & P_{n+1}T\ECL{\pinC{n}X}{n+2}
}
$$
\noindent is a homotopy pullback square in \ $\D/T\BCL$, \ where \ 
$\Lambda:=\pinC{0}X$. 
\end{lemma}

\begin{proof}
Set \ $E:=T\ECL{J}{n+1}$, \ $M^{n-1}:=TP_{n-1}X$, \ and \ $M^{n}:=TP_{n}X$. \ 
The naturality of the comparison exact sequence, applied to the maps
in \eqref{ezero}, \ combined with Fact \ref{fsix}, imply that the
vertical maps in the following commutative diagram are isomorphisms:
$$
\xymatrix{
\pinD{n+2}E \rto \dto^{\cong} & \pinD{n+1}M^{n} \rto \dto^{\cong} &
\pinD{n+1}M^{n-1} \rto^>>>>{Tk_{n-1}} \dto^{\cong} & 
\pinD{n+1}E \rto \dto^{\cong} &
\pinD{n}M^{n} \rto \dto^{\cong} & \pinD{n}M^{n-1} \dto^{\cong} \\
0\rto & \Cok h^{X}_{n+1} \ar@{|->}[r] & \Gamma_{n}X \rto^{s^{X}_{n}} &
\pinC{n}X \rto^{h^{X}_{n}} & \pinD{n}T X \ar@{->>}[r] & \Cok h_{n}^{T}
}
$$
\noindent and since the bottom row is part of the comparison long
exact sequence, and the rest of the top sequence to the right 
is exact for by \eqref{eeight}, the $k$-invariant square \eqref{ezero} 
induces a long exact sequence after applying \ $\pin$ \ (except in the
bottom dimensions). \ The obvious map from \ $M^{n}$ \ to the fiber of \
$Tk_{n-1}$ \ is thus a weak equivalence in \ $\D/T\BCL$ \ through
dimension \ $n+1$.
\end{proof}

%
%
\begin{cor}\label{czero}\stepcounter{subsection}
For \ $T:\C\to\D$ \ as above, \ for any \ $X\in\C$ \ and \ $n\geq 1$ \ 
the natural map \ $r\q{n}:X\to P_{n}X$ \ of \S \ref{sps} induces an \ 
isomorphism \ $\Gamma_{k}X\cong\Gamma_{k}P_{n}X$ \ for \ $k\leq n+1$.
\end{cor}

\begin{proof}
For each \ $A\in\A$, \ take fibers vertically and horizontally of
the commutative square:
$$
\xymatrix{
\map_{\BCL}(A,P_{n}X) \rto^-{\protect{h^{P_{n}X}}} \dto^{(k_{n})_{\ast}} & 
\map_{\BDL}(TA,TP_{n}X) \dto^{(Tk_{n})_{\ast}} \\ 
\map_{\BCL}(A,\ECL{\pinC{n+1}X}{n+2}) \rto^-{\protect{h^{E}}} & 
\map_{\BDL}(TA,T\ECL{\pinC{n+1}X}{n+2})~,\\
}
$$
and use Lemma \ref{lthree} and \S \ref{dspecial}(iv).
\end{proof}

\begin{remark}\label{rhem}\stepcounter{subsection}
For \ $\C=s\G$ \ this follows from the fact that \
$\Gamma_{n}X\cong\Omega\pin_{n-1}X$, \ while for the algebraic cases
of \S \ref{ermc}(i-ii), this follows from the fact that \
$H_{n+1}(K(\pi,n);\bZ)=0$ \ for \ $n\geq 1$.
\end{remark}

\subsection{The extension}
\label{sext}\stepcounter{thm}

The map \ $r\q{n}:X\to P_{n}X$ \ induces a map of comparison exact
sequences: 
%
\setcounter{equation}{\value{thm}}\stepcounter{subsection}
\begin{equation}\label{etwentysix}
\xymatrix{
\pinC{n+2}X \dto \rto^{h^{X}_{n+2}} & 
\pinC{n+2}TX \dto^{\pinD{n+2}Tr\q{n}} \rto^{\partial^{\star}_{n+2}} & 
\Gamma_{n+1}X \dto^{=} \rto^{s_{n+1}} & \pinC{n+1}X \dto \rto^{h^{X}_{n+1}} & 
\pinD{n+1}TX \dto^{\pinD{n+1}Tr\q{n}} \rto^{\partial^{\star}_{n+1}} & 
\Gamma_{n}X \dto^{=} \\
0 \rto & \pinD{n+2}M^{n} \rto^{\cong} & \Gamma_{n+1}P_{n}X \rto & 0 \rto &
\pinD{n+1}M^{n} \rto & \Gamma_{n}P_{n}X
}
\end{equation}
\setcounter{thm}{\value{equation}}
\noindent so that \ $\pinC{n+1}X$ \ fits into a short exact sequence 
of \PAa s:
%
\setcounter{equation}{\value{thm}}\stepcounter{subsection}
\begin{equation}\label{eten}
0\to \Cok\pinD{n+2}Tr\q{n} \to \pinC{n+1}X \to
\Ker\pinD{n+1}Tr\q{n} \to 0,
\end{equation}
\setcounter{thm}{\value{equation}}
\noindent where 
%
\setcounter{equation}{\value{thm}}\stepcounter{subsection}
\begin{equation}\label{etwentyone}
\Cok\pinD{n+2}Tr\q{n}\cong \Ker h^{X}_{n+1}\hsm \text{and}\hsm 
\Ker\pinD{n+1}Tr\q{n} \cong \Image h^{X}_{n+1}.
\end{equation}
\setcounter{thm}{\value{equation}}

Since \ $h^{X}_{n+1}$ \ is a map of modules over \ $\Lambda:=\pinC{0}X$, \ 
by Theorem \ref{tone}, \eqref{eten} is actually a short exact sequence
of modules over $\Lambda$, and we can classify the possible values of \ 
$J\in\RM{\Lambda}$ \ (the candidates for \ $\pinC{n+1}X$) \ using the following:
%
%
\begin{prop}\label{peleven}\stepcounter{subsection}
Given \ $Tr\q{n}:T X\to T P_{n}X$, \ a choice for the 
isomorphism class of \ $\pinC{n+1}X$ \ uniquely determines an element of \ 
$$
\Ext_{\RM{\Lambda}}(\Ker(Tr\q{n})_{n+1}, \Cok(Tr\q{n})_{n+2}).
$$
\end{prop}

\begin{proof}
Since \ $\RM{\Lambda}$ \ is an abelian category, with a set  \ 
$\{A_{\ab}\otimes S^{n}\amalg\BDL\}_{A\in\A,n\in\bN}$ \ 
of projective generators, the argument of \cite[III]{MacH} carries
over to our setting. 
\end{proof}

\begin{remark}\label{rinfo}\stepcounter{subsection}
Observe that given \ $P_{n}X$, \ we know the comparison exact sequence 
\eqref{efive} for \ $X$ \ only from \ 
$s_{n}:\Gamma_{n-1}X\to\pinC{n}X$ \ down. However, if \ 
$\pinD{i}Tr\q{n}:\pinD{i}TX\to\pinD{i}M^{n}$ \ (for \ $i\geq 0$) \ and the 
extension \eqref{eten} are also known,  all we need in order to
determine \eqref{efive} for $X$ from \ 
$\partial_{n+3}^{\star}:\pinD{n+3}TX\to\Gamma_{n+1}X$ \ down is the 
homomorphism
$$
\pinD{n+3}Tr\q{n+1}:\pinD{n+3}TX\to\pinD{n+3}TP_{n+1}X~,
$$
which is just \ $\partial_{n+3}^{\star}$, \ as one can see 
from \eqref{etwentysix}.
\end{remark}

%
%
\begin{prop}\label{ptwelve}\stepcounter{subsection}
For any \ $\Lambda\in\D$, \ $J',J''\in\RM{\Lambda}$, \ and \ 
$n\geq 2$, \ there is a natural isomorphism \ 
$\Ext_{\RM{\Lambda}}(J'',J')\cong\HL{n+1}{\EDL{J''}{n}}{J'}$.
\end{prop}

In particular, this implies that \ $\HL{n+1}{\EDL{-}{n}}{-}$ \ is 
\emph{stable} \ -- \ i.e., independent of $n$.

\begin{proof}
By Proposition \ref{ptwo}\emph{ff.} there is a natural isomorphism \ 
$$
\HL{n+1}{\EDL{J''}{n}}{J'}
\cong[\EDL{J''}{n},\EDL{J'}{n+1}]_{s\D/\BDL},
$$
\noindent and given a map \ $\psi:\EDL{J''}{n}\to\EDL{J'}{n+1}$, \ we can 
form the fibration sequence over \ $\BDL$ \ (that is, pullback square
as in \eqref{ezero}):
$$
\Omega\EDL{J''}{n}\xra{\Omega\psi}\Omega\EDL{J'}{n+1}\simeq\EDL{J'}{n}
\to F\to\EDL{J''}{n}\xra{\psi}\EDL{J'}{n+1}.
$$
\noindent From the corresponding long exact sequence in homotopy 
for this sequence in $\D$, we obtain a short exact sequence of 
modules over $\Lambda$:
%
\setcounter{equation}{\value{thm}}\stepcounter{subsection}
\begin{equation}\label{eeleven}
0\to J'\to J\to J''\to 0.
\end{equation}
\setcounter{thm}{\value{equation}}

On the other hand, given a short exact sequence \eqref{eeleven} in \ 
$\RM{\Lambda}$, \ we can construct a map \ 
$\psi:\EDL{J''}{n}\to\EDL{J'}{n+1}$ \ over \ $\BDL$ \ as follows:

Assume \ $E:=\EDL{J''}{n}$ \ is constructed  starting with \ 
$\sk{n-1}\EDL{J''}:=\sk{n-1}\BDL$, \ and \ 
$E_{n}\simeq W\amalg L_{n}\BDL$ \ (cf.\ \S \ref{dlat}), \ where $W$ is free, 
equipped with a surjection \ $\phi:W\to J''$. \ Because \ $J\epic J''$ \ is 
a surjection, and $W$ is free, we can lift $\phi$ to \ 
$\phi':W\to J$, \ defining a map \ $\tilde{\phi}':Z_{n}\EDL{J''}{n}\to J$. \ 
Since \  $\pinD{n}\EDL{J''}{n}=J''$, \ the restriction of \ 
$\tilde{\phi}'$ \ to \ $B_{n}\EDL{J''}{n}=\Ker\{Z_{n}\EDL{J''}{n}\to J''\}$ \ 
factors through \ $\psi:B_{n}\EDL{J''}{n}\to J'=\Ker\{J\to J''\}$. \
Precomposing with \ $\bd:C_{n+1}\EDL{J''}{n}\to B_{n}\EDL{J''}{n}$ \ defines \ 
$\psi:\EDL{J''}{n}\to\EDL{J'}{n+1}$, \ which classifies \eqref{eeleven} 
as before.
\end{proof}

%
%
\begin{cor}\label{cone}\stepcounter{subsection}
For $\Lambda$, $J'$, and $J''$ as above, there is a natural isomorphism:
$$
\Ext_{\RM{\Lambda}}(J'',J')\cong\HL{n+1}{\ECL{J''}{n}}{J'}.
$$
\end{cor}

\begin{proof}
This follows from \eqref{etwentyfour}-\eqref{enine} and the 
naturality of \ $P^{\D}_{n+1}$.
\end{proof}

\begin{defn}\label{dmps}\stepcounter{subsection}
Given \ $X\in\C$, \ its \emph{$n$-th modified Postnikov section}, 
denoted by \ $\bPa{n}X$, \ is defined as follows\vsm :

Let \ $\K:=\{f:A\otimes S^{n+1}\to X~|\ 
A\in\A,~[f]\in\Ker h_{n+1}^{T}\subset\pinC{n+1}X\}$, \ and let
$C$ be the cofiber of the obvious map \ 
$\Phi:\bigvee_{f\in\K}~A\otimes S^{n+1}\to X$ \  
(so that \ $\pinC{n+1}C\cong\Cok\Phi$), \ with \ $\bPa{n}X:= P_{n+1}C$. \ 
There are then natural maps \ $\bpa{n+1}:P_{n+1}X\to\bPa{n}X$ \ 
(induced by \ $X\to C$), \ as well
as \ $\bpc{n}:\bPa{n}X\to P_{n}X$ \ (which is just \ 
$p\q{n}_{C}:P_{n+1}C\to P_{n}C\cong P_{n}X$), \ with \ 
$\bpc{n}\circ\bpa{n}=p\q{n}_{X}:P_{n+1}X\to P_{n}X$. \ 
Note that \ $\pinC{n+1}\bPa{n}X\cong\Image h^{X}_{n+1}$, \ and \ 
$P_{n}\bPa{n}X\cong P_{n}X$.
\end{defn}

The map \ $\brp{n}:=\bpa{n}\circ r\q{n}:X\to\bPa{n}X$ \ induces a map of 
comparison exact sequences:
$$
\xymatrix{
\pinC{n+2}X \rto^{h^{X}_{n+2}} \dto & 
\pinD{n+2}T X \rto^{\partial^{\star}_{n+2}} \dto^{\pinD{n+2}T\brp{n}} &
\Gamma_{n+1}X \rto^{s_{n+1}} \dto^{=} &
\pinC{n+1}X \rto^{h^{X}_{n+1}} \dto &
\pinD{n+1}T X \rto^{\partial_{n+1}} \dto_{=}^{\pinD{n+1}T\brp{n}}
& \Gamma_{n}X \dto^{=} \\
0 \rto & \pinD{n+2}T\bPa{n}X \rto^{\cong} & \Gamma_{n+1}\bPa{n}X \rto^{0} & 
\pinC{n+1}\bPa{n}X \ar@{|->}[r] & \pinD{n+1}T\bPa{n}X \rto & \Gamma_{n}\bPa{n}X
}
$$
\noindent so that:
%
\setcounter{equation}{\value{thm}}\stepcounter{subsection}
\begin{equation}\label{etwenty}
\pinD{k}T\bPa{n}X\cong\begin{cases}
        \pinD{k}TX& \text{for \ } k\leq n+1,\\
        \Gamma_{n+1}X & \text{for \ }k=n+2,\\
        \Gamma_{k-1}\bPa{n}X & \text{for \ }k\geq n+3~.
\end{cases}
\end{equation}
\setcounter{thm}{\value{equation}}

Thus \ $\brp{n}$ \ induces a weak equivalence \ $P_{n+1}TX\simeq
P_{n+1}T\bPa{n}X$, \ which, together with the existence of the
appropriate maps \ $P_{n+1}X\xra{\bpa{n}}\bPa{n}X\xra{\bpc{n}}P_{n}X$, \ 
determines \ $\bPa{n}X$ \ up to homotopy. In fact we have:

%
%
\begin{prop}\label{pfive}\stepcounter{subsection}
$\bPa{n}X$ \ is determined uniquely (up to weak equivalence) by \ 
$P_{n}X$ \ and the map \ $\rho:=P_{n+1}Tr\q{n}:P_{n+1}TX\to P_{n+1}TP_{n}X$.
\end{prop}

\begin{proof}
Note that \ $I_{n+1}:=\Ker \pinD{n+1}\rho$ \ is isomorphic to \ 
$\Image h^{X}_{n+1}$ \ and \ $C_{n+1}:=\Image \pinD{n+1}\rho$ \ 
is isomorphic to \ $\Cok h^{X}_{n+1}$ \ by \eqref{etwentyone}. \ 

We construct \ $Y\simeq\bPa{n}X$ \ as follows, starting with \ 
$\sk{n+1}Y:=\sk{n+1}P_{n}X$; \ by Remark \ref{rcsk}, we may 
assume \ $\sk{n+1}TX=\sk{n+1}TP_{n}X$, \ so that \ 
$P_{n}TX\cong P_{n}TP_{n}X$. \ By Fact \ref{ffive}), the lower right 
hand square in Figure \ref{fig6} commutes in $\D$, thus inducing the
rest of the diagram, in which the rows and columns are fibration
sequences over \ $\BDL$. 

\begin{center}
%
%
\setcounter{figure}{\value{thm}}\stepcounter{subsection}
\begin{figure}[hbtp]
\begin{picture}(300,110)(25,-5)
%
%
\put(15,100){$F$}
\multiput(35,105)(3,0){34}{\circle*{.5}}
\put(137,105){\vector(1,0){2}}
\put(85,110){$\hat{\rho}$}
\put(145,100){$P_{n+1}TP_{n}X$}
\multiput(205,105)(3,0){13}{\circle*{.5}}
\put(245,105){\vector(1,0){2}}
\put(217,110){$\tkp{n}$}
\put(255,100){$\EDL{I_{n+1}}{n+2}$}
%
%
\put(20,92){\vector(0,-1){30}}
\put(8,75){$\simeq$}
\put(25,75){$\lambda$}
\put(170,92){\vector(0,-1){30}}
\put(175,75){$p\q{n}_{TP_{n}X}$}
\put(290,92){\vector(0,-1){30}}
\put(295,75){$i_{\ast}$}
%
%
\put(0,50){$P_{n+1}TX$}
\put(45,55){\vector(1,0){55}}
\put(65,62){$p\q{n}_{TX}$}
\put(105,50){$P_{n}TX\cong P_{n}TP_{n}X$}
\put(205,55){\vector(1,0){45}}
\put(220,61){$k_{n}^{TX}$}
\put(255,50){$\EDL{\pinD{n+1}TX}{n+2}$}
%
%
\put(20,42){\vector(0,-1){30}}
\put(170,42){\vector(0,-1){30}}
\put(175,28){$k_{n}^{TP_{n}X}$}
\put(290,42){\vector(0,-1){30}}
\put(295,28){$q_{\ast}$}
%
%
\put(10,0){$\BDL$}
\put(55,5){\vector(1,0){65}}
\put(125,0){$\EDL{C_{n+1}}{n+2}$}
\put(215,5){\vector(1,0){35}}
\put(230,10){$=$}
\put(255,0){$\EDL{C_{n+1}}{n+2}$}
\end{picture}
\caption[fig6]{}
\label{fig6}
\end{figure}
\setcounter{thm}{\value{figure}}
\end{center}

In particular, the induced map \ 
$\tkp{n}:P_{n+1}TP_{n}X\to\EL{I_{n+1}}{n+2}$ \ provides a canonical lifting of:
$$
k_{n}^{TX}\circ p\q{n}_{TP_{n}X}:P_{n+1}TP_{n}X\to\EDL{\pinD{n+1}TX}{n+2}
$$
\noindent to \ $\EDL{I_{n+1}}{n+2}$. \ Composing it with the natural map \ 
$r\q{n+1}:TP_{n}X\to P_{n+1}TP_{n}X$ \ defines an element in: 
$$
[TP_{n}X,\EDL{I_{n+1}}{n+2}]\cong\HL{n+2}{P_{n}X}{I_{n+1}}~,
$$
\noindent which we call the \emph{$n$-th modified $k$-invariant} for $X$\vsm.

If \ $\hk{n}:P_{n}X\to\ECL{I_{n+1}}{n+2}$ \ is the map corresponding to
$\tkp{n}$ \ under \eqref{etwentyfour}), then its homotopy fiber $Y$ is 
(weakly equivalent to) \ $\bPa{n}X$, \ as one can verify by calculating \ 
$\pinC{\ast}Y$. \ Note that Lemma \ref{lthree} implies that \ 
$F\simeq P_{n+1}TP_{n}X$, \ so that $\lambda$ is the homotopy inverse
of the weak equivalence \ $P_{n+1}\rho:TX\to TP_{n}X$, \ which
completes the construction. 
\end{proof}

\begin{remark}\label{rmin}\stepcounter{subsection}
Note that there is a certain indeterminacy in our description of \ 
$\tkp{n}$, \ and thus of \ $\hk{n}$, \ since we must make the lower right
corner of Figure \ref{fig6} into a strict commuting diagram of fibrations,
rather than one which commutes only up to homotopy. However,
\end{remark}

%
%
\begin{fact}\label{fseven}\stepcounter{subsection}
The indeterminacy for \ $\tkp{n}$ \ as an induced map is contained in 
the indeterminacy for \ $\tkp{n}$ \ as a $k$-invariant
for \ $P_{n+1}TX=P_{n+1}TY$.
\end{fact}

\begin{proof}
Let \ $M:=TP_{n}X$. \ Making the lower right corner of Figure
\ref{fig6} commute on the nose (assuming \ $q_{\ast}$ \ is already a
fibration) requires the choice of a homotopy
$$
H:P_{n}TX\to\Omega\EDL{C_{n+1}}{n+2}=\EDL{C_{n+1}}{n+1}~,
$$
\noindent so the indeterminacy for \ $\tkp{n}$ \ as defined above is \ 
$\psi_{\ast}p^{\ast}[P_{n}TX,\EDL{C_{n+1}}{n+1}]$, \ where \ 
$\psi:\EDL{C_{n+1}}{n+1}\to\EDL{I_{n+1}}{n+2}$ \ classifies the extension \ 
$$
0\to I_{n+1}\to\pinD{n+1}TX\to C_{n+1}\to 0
$$
\noindent (Proposition \ref{ptwelve}), and \ 
$p=p\q{n}_{M}:P_{n+1}M\to P_{n}M=P_{n}TX$. 

On the other hand, the $k$-invariant \ 
$\tkp{n}^{M}:P_{n+1}M\to\EDL{I_{n+1}}{n+2}$ \ 
for \ $P_{n+1}TP_{n}X$ \ (which is \ $P_{n+1}TX$) \ is 
determined only up to the actions of the group \ $\haut_{\Lambda}(P_{n+1}M)$ \
of homotopy self-equivalences of \ $P_{n+1}M$ \ over \
$\BDL$, \ and of \ $\AutL(I_{n+1})$, \ the group of automorphisms of
modules over $\Lambda$ of \ $I_{n+1}$, \ in \
$[P_{n+1}M,\EDL{I_{n+1}}{n+2}]$. 
Thus given a map \ $f:P_{n}M\to\EDL{C_{n+1}}{n+1}$, \ we obtain a 
self-map \ $g:P_{n+1}M\to P_{n+1}M$ \ such that \ 
$P_{n}g=\Id_{P_{n}M}$ \ and \ $\pinD{n+1}g=\Id$, \ by letting \ 
$g=\Id+i_{\ast}p^{\ast}(f)$, \ for \ $i:\EDL{C_{n+1}}{n+1}\to P_{n+1}M$ \ 
the inclusion of the fiber. It is readily verified that 
$g$ induces the identity on \ $\pinD{\ast}P_{n+1}M$, \ so \ 
$[g]\in \haut_{\Lambda}(P_{n+1}M)$, \ and that \ 
$\tkp{n}+\psi_{\ast}p^{\ast}(f)$ \ is obtained from \ $\tkp{n}$ \ under
the action of \ $[g]$ \ on \ $\HL{n+2}{P_{n+1}M}{I_{n+1}}$.
\end{proof}

\begin{notation}\label{nbpa}\stepcounter{subsection}
Given \ $W\simeq P_{n}X$ \ and \ $\rho:P_{n+1}TX\to P_{n+1}TW$, \ 
Proposition \ref{pfive} allows us to write \ $\bPa{n}(W,\rho)$, \ or
simply \ $\bPa{n}W$ \ for \ $\bPa{n}X\in\C$, \ which
they determine up to homotopy. This comes equipped with a weak equivalence \ 
$\rho:P_{n+1}TX\to P_{n+1}T\bPa{n}W$ \ lifting $\rho$.
\end{notation}

%
%
\begin{cor}\label{ctwo}\stepcounter{subsection}
The weak equivalence \ $\rho:P_{n+1}TX\to P_{n+1}T\bPa{n}W$ \ is
well-defined up to homotopy.
\end{cor}

\begin{proof}
The map $\rho$ is inverse to $\lambda$ in Figure \ref{fig6}, which is induced
by the upper right hand square, which is determined by \ $\tkp{n}$ \ 
and thus up to a self-equivalence \ $g:P_{n+1}TW\to P_{n+1}TW$, \ 
according to Fact \ref{fseven}. \ But such a $g$ induces a canonical
self-equivalence \ $g':F'\to F$, \ where \ 
$F':=\Fib(\tkp{n}\circ g)$, \ and the resulting \ 
$\lambda':F'\simeq P_{n+1}TX$ \ satisfies \ 
$\lambda\circ g'\simeq\lambda'$.
\end{proof}

\begin{defn}\label{dae}\stepcounter{subsection}
For $W\simeq P_{n}X$ \ and \ $\rho:P_{n+1}TX\to P_{n+1}TW$ \ as above, 
an extension
%
\setcounter{equation}{\value{thm}}\stepcounter{subsection}
\begin{equation}\label{ethirteen}
0\to \Cok \pinD{n+2}\rho\hra J\epic\Ker\pinD{n+1}\rho \to 0
\end{equation}
\setcounter{thm}{\value{equation}}
\noindent is called \emph{allowable} if its classifying cohomology class \ 
$$
[\psi]\in\HL{n+3}{\EDL{\Cok \pi_{n+2}\rho}{n+2}}{\Ker \pi_{n+1}\rho }
$$
(cf.\ Proposition \ref{ptwelve}) satisfies \ $[\psi]\circ\hk{n}=0$.
\end{defn}

%
%
\begin{prop}\label{pthirteen}\stepcounter{subsection}
For any \ $X\in\C$, \ the extension \eqref{eten} is
allowable.
\end{prop}

\begin{proof}
Writing \ $V\simeq P_{n+1}X$ \ and \ $Y\simeq\bPa{n}X$, \ by
naturality we have a commutative square:
$$
\xymatrix{
P_{n}V \rto^-{\protect{k_{n}}} \dto^{=} & \ECL{\pinC{n+1}V}{n+2}
\dto^{q_{\ast}} \\ 
P_{n}Y \rto^-{\protect{k_{n}}} & \ECL{\Ker\pi_{n+1}Tr\q{n}}{n+2}. \\
}
$$

Lemma \ref{lthree} and \eqref{etwentyfour} then yield the following
commuting diagram in $\D$ in which the rows and columns are all
fibration sequences over \ $\BDL$: 

\begin{center}
%
%
\begin{picture}(380,170)(0,-50)
%
%
\put(0,100){$\EDL{\Ker h_{n}}{n+1}$}
\put(100,105){\vector(1,0){45}}
\put(150,100){$\BDL$}
\put(180,105){\vector(1,0){40}}
\put(225,100){$\EDL{\Ker h_{n+1}}{n+2}$}
%
%
\put(50,95){\vector(0,-1){33}}
\put(160,95){\vector(0,-1){30}}
\put(270,95){\vector(0,-1){30}}
%
%
\put(30,50){$TP_{n+1}X$}
\put(79,55){\vector(1,0){60}}
\put(95,62){$Tr\q{n}$}
\put(143,50){$TP_{n}X$}
\put(182,55){\vector(1,0){40}}
\put(195,61){$k$}
\put(225,50){$\EDL{\pinC{n+1}X}{n+2}$}
%
%
\put(50,45){\vector(0,-1){33}}
\put(55,28){$r\q{n+2}$}
\put(160,45){\vector(0,-1){30}}
\put(165,28){$=$}
\put(270,45){\vector(0,-1){30}}
\put(255,28){$q_{\ast}$}
%
%
\put(40,0){$TY$}
\put(70,5){\vector(1,0){70}}
\put(143,0){$TP_{n}X$}
\put(177,5){\vector(1,0){40}}
\put(200,11){$\tkp{}$}
\put(225,0){$\EDL{\Image h_{n+1}}{n+2}$}
\put(270,-7){\vector(0,-1){23}}
\put(257,-18){$\psi$}
\put(225,-43){$\EDL{\Ker h_{n+1}}{n+3}$}
\end{picture}
\end{center}

\noindent The map $k$ is induced by \ $k_{n}$, \ and \ $\tkp{}$ \ is
induced by \ $\hk{n}$. \ The claim then follows from the
universal property for fibrations. 
\end{proof}

%
%
\sect{The fiber of a special spherical functor}
\label{cfib}

Let  \ $T:\C\to\D$ \ be a special spherical functor.
We would like to use the results of Section \ref{ccps}  in order to
determine whether a given \ $G\in\D$ \ is (up to homotopy) of the form \ 
$TX$ \ for some \ $X\in\C$ \ -- \  and if so, how we can distinguish 
between such \emph{realizations}, or liftings. 

\subsection{Lifting objects of $\D$}
\label{srsp}\stepcounter{thm}

Let us assume for simplicity that \ $\Lambda:=\pinD{0}G$ \ is a \PCa,
and that the map \ $\phi_{T}:\Lambda\to\pinD{0}T\BCL$ \ of \S
\ref{dspecial}(i) is an isomorphism. [In the general case, we are
faced with an additional, purely algebraic, problem of determining the
fiber of the functor \ $T_{\ast}:\PCAlg\to\PDAlg$ \ 
(compare \cite{BPescF}); \ we bypassed this question in \S \ref{dsfunc}(iv).  
   
We want a map \ $\varphi:TX\to G$ \ inducing isomorphisms \
$\pinD{i}TX\to\pinD{i}G$ \ for \ $i\geq 0$. \ 
Our approach is inductive: we are trying to define a tower in $\C$:
%
\setcounter{equation}{\value{thm}}\stepcounter{subsection}
\begin{equation}\label{eseven}
\dotsb \xra{p\q{n+1}} \Xpn{n+1} \xra{p\q{n}} \Xpn{n} \xra{p\q{n-1}} \dotsb
\xra{p\q{0}} \Xpn{0}\simeq\BCL
\end{equation}
\setcounter{thm}{\value{equation}}
\noindent which are to serve as the modified Postnikov tower of the
(putative) \ $X\in\C$ -- \ so that in the end we will have \ 
$X:=\holim_{n}\Xpn{n}$.

At the $n$-th stage \ ($n\geq 0$), we have constructed \ $\Xpn{n}$ \ as our 
candidate for \ $\bPa{n}X$ \ -- \ so in particular if we let \ 
$\Xn{n}:=P_{n}\Xpn{n}$, \ (our candidate for the ordinary $n$-th Postnikov 
section of \ $X$), \ then \ $T\Xn{n}$ \ satisfies \eqref{enine}, \ 
$T\Xpn{n}$ \ satisfies \eqref{etwenty}, \ and of course \ 
$\Xpn{n}=P_{n+1}\Xpn{n}$.

Assume also, as part of our inductive hypothesis, a given weak equivalence:
%
\setcounter{equation}{\value{thm}}\stepcounter{subsection}
\begin{equation}\label{eseventeen}
\hr{n}:P_{n+1}G\xra{\simeq}P_{n+1}T\Xpn{n}.
\end{equation}
\setcounter{thm}{\value{equation}}

We start the induction with \ $\Xn{0}:=\BCL$. \ The natural
map \ $r\q{1}:G\to P_{1} T\BCL=\BDL$ \ allows us to define \
$\Xpn{0}$, \ together with \ $\hr{0}:P_{1}G\xra{\simeq}P_{1}T\Xpn{0}$, \ 
as in Definition \ref{dmps} (see \S \ref{rmin}).

\subsection{Lifting $\rho\q{n}$}
\label{slrpn}\stepcounter{thm}

The first stage in the inductive step occurs in $\D$: we must lift \ 
$\hr{n}$ \ to \ $\rho\q{n}:P_{n+2}G\to P_{n+2}T\Xpn{n}$. \ 
Note that by Remark \ref{rinfo} and Fact \ref{fsix}, we already know the 
comparison exact sequence \eqref{efive} for the putative $X$ from \ 
$h_{n+1}$ \ down; the lifting \ $\rho:=\rho\q{n}$ \ will determine \ 
$\partial_{n+2}:\pinD{n+2}G\to\Gamma_{n+1}\Xpn{n}$ \ in addition, 
since this is just \ $\pi_{n+2}\rho$, \ so that \ 
$C_{n+2}:=\Image \pinD{n+2}\rho$ \ is our candidate for \ 
$\Cok h_{n+2}^{X}$, \ while \ $K_{n+1}:=\Cok \pinD{n+2}\rho$ \ is our 
candidate for \ $\Ker h_{n+1}^{X}$.

From \eqref{etwenty} we see that the obstruction is the class:
%
\setcounter{equation}{\value{thm}}\stepcounter{subsection}
\begin{equation}\label{enineteen}
\chi_{n}:=k_{n+1}^{T\Xpn{n}}\circ\rho\q{n}\in
\HL{n+3}{G}{\Gamma_{n+1}\Xpn{n}}~,
\end{equation}
\setcounter{thm}{\value{equation}}
\noindent and the different liftings are classified by \ 
$\HL{n+2}{G}{\Gamma_{n+1}\Xpn{n}}$.

\subsection{Constructing \ $\Xn{n+1}$}
\label{scpx}\stepcounter{thm}

The next step is to choose a cohomology class \ 
$\hk{n}$ \ in \ $\HL{n+2}{\Xpn{n}}{K_{n+1}}$. \ This fits into a
commutative diagram with rows and fibers all fibration sequences over \ $\BCL$:
$$
\xymatrix{
\BCL \rto \dto & \ECL{I_{n+1}}{n+1} \rto^{=} \dto^{i} &
\ECL{I_{n+1}}{n+1} \dto^{\psi} \\
\Xn{n+1} \rto^{\bpa{n}} \dto & \Xpn{n} \rto^>>>>>>>{\hk{n}}\dto & 
\ECL{K_{n+1}}{n+2} \dto^{j_{\ast}} \\
\Xn{n+1}\rto & \Xn{n} \ar@{.>}[r]^>>>>>>>>>>{k_{n}} & \ECL{J }{n+2}
}
$$
for the bottom fibration sequence \ $\Xn{n+1}\to\Xn{n}\to\ECL{J }{n+2}$ \ 
as indicated (though we shall not need this).

Note that $J$, our candidate for \ $\pinC{n+1}X$, \ fits into the short 
exact sequence of modules over $\Lambda$:
$$
0\to K_{n+1}\hra J\epic I_{n+1}\to 0,
$$
\noindent as in \eqref{eten}, and is classified by \ 
$\psi:=\hk{n}\circ i\in\HL{n+2}{\ECL{I_{n+1}}{n+1}}{K_{n+1}}$, \ 
as in Corollary \ref{cone}. Moreover, this extension is
obviously allowable in the sense of \S \ref{dae}.

\subsection{Lifting $\rho$}
\label{slr}\stepcounter{thm}

To complete the induction on \eqref{eseventeen}, we must lift \ 
$\rho:G\to P_{n+2}T\Xpn{n}$. \ This will be done in two steps\vsm :

First, note that we obtain a commuting diagram:

\begin{center}
%
%
$$
\begin{picture}(350,180)(0,0)
%
%
\put(15,170){$P_{n+2}G$}
\put(60,175){\vector(1,0){205}}
\put(140,180){$\rho$}
\put(270,170){$P_{n+2}T\Xpn{n}$}
\multiput(60,168)(2,-1){25}{\circle*{.5}}
\put(110,143){\vector(2,-1){2}}
\put(78,143){$\rho$}
\put(115,135){$P_{n+2}T\Xn{n+1}$}
\put(205,142){\vector(2,1){55}}
\put(230,143){$\tilde{i}_{\ast}$}
%
%
\put(35,162){\vector(0,-1){62}}
\put(3,130){$p\q{n+1}_{G}$}
\put(150,128){\vector(0,-1){30}}
\put(155,112){$p\q{n+1}_{T\Xn{n+1}}$}
\put(290,162){\vector(0,-1){62}}
\put(295,130){$p\q{n+1}_{T\Xpn{n}}$}
%
%
\put(15,85){$P_{n+1}G$}
\put(60,90){\vector(1,0){50}}
\put(80,95){$f$}
\put(80,82){$\simeq$}
\put(115,85){$P_{n+1}T\Xn{n+1}$}
\put(208,90){\vector(1,0){55}}
\put(230,95){$g$}
\put(230,82){$\simeq$}
\put(270,85){$P_{n+1}T\Xpn{n}$}
%
%
\put(35,77){\vector(0,-1){62}}
\put(15,40){$k_{n+1}^{G}$}
\put(150,77){\vector(0,-1){28}}
\put(155,60){$k_{n+1}^{T\Xn{n+1}}$}
\put(290,77){\vector(0,-1){62}}
\put(295,40){$k_{n+1}^{T\Xpn{n}}$}
%
%
\put(-30,0){$\EDL{\pin_{n+2}G}{n+3}$}
\put(75,5){\vector(1,0){190}}
\put(145,11){$(\pin_{n+2}\rho)_{\ast}$}
\multiput(70,10)(2,1){23}{\circle*{.5}}
\put(116,33){\vector(2,1){2}}
\put(90,30){$q_{\ast}$}
\put(123,35){$\EDL{C_{n+2}}{n+3}$}
\put(213,37){\vector(2,-1){53}}
\put(235,30){$i_{\ast}$}
\put(270,0){$\EDL{\pinD{n+1}TX}{n+3}$}
\end{picture}
$$
%
%
\end{center}

\noindent in which the columns are fibration sequences over \
$\BCL$, \ since by definition \ 
$$
\pinD{n+2}\rho:\pinD{n+2}G \to\pinD{n+1}T\Xpn{n}=\pinD{n+1}TX
$$
factors through \ $C_{n+2}:=\Image \pinD{n+2}\rho$, \ so that the
bottom triangle  commutes. 

Since the natural $K$-invariant \ \ $k_{n+1}^{G}$ \ is given,
the other two $k$-invariants in the diagram above are determined by
inverting the given homotopy equivalences \ 
$f:P_{n+1}G\to P_{n+1}T\Xn{n+1}$ \ and \ 
$g:P_{n+1}G\to P_{n+1}T\Xpn{n}$ \ (assuming all objects in $\D$ are
fibrant and cofibrant), \ and letting \ 
$k_{n+1}^{T\Xn{n+1}}:=q_{\ast}\circ k_{n+1}^{G}\circ f^{-1}$ \ and \ 
$k_{n+1}^{T\Xpn{n}}:=i_{\ast}\circ k_{n+1}^{G}\circ g^{-1}$, \ 
using Fact \ref{ffive}.

Therefore, the map \ $\rho:G\to P_{n+2}T\Xpn{n}$ \ 
lifts to \ $\rho:P_{n+2}G\to P_{n+2}T\Xn{n+1}$ \ (which is induced by \ 
$q_{\ast}$). \ In fact, the lifting $\rho$ is unique up to homotopy.
Moreover, from the proof of Proposition \ref{pfive} we see that this suffices 
to define \ $\Xpn{n+1}$, \ as well as determining a lifting of $\rho$ to 
a weak equivalence \ $\hr{n+1}:P_{n+2}G\to P_{n+2}T\Xpn{n+1}$\vsm . 

We may summarize our results in:

%
%
\begin{thm}\label{tfour}\stepcounter{subsection}
Given \ $G\in\D$, \ there is an object \ $X\in\C$ \ 
such that \ $TX\simeq G$ \ if and only if there is a tower as in 
\eqref{eseven}, serving as the modified Postnikov tower for $X$. \ If we 
have constructed \ $\Xpn{n}$ \ satisfying \eqref{eseventeen} for \ $n$, \
a necessary and sufficient condition for the existence of an \ $\Xpn{n+1}$ \ 
satisfying \eqref{eseventeen} for \ $n+1$ \ is the vanishing of \ 
$\chi_{n}\in\HL{n+3}{G}{\Gamma_{n+1}\Xpn{n}}$. \ The choices are classified by: 
\begin{enumerate}
\renewcommand{\labelenumi}{$\bullet$~}
\item $\HL{n+2}{G}{\Gamma_{n+1}{\Xpn{n}}}$ \ (distinguishing the
  liftings of \ $\hr{n}$ \ to \ $P_{n+2}T\Xpn{n}$); \ and 
\item $\hk{n}\in\HL{n+2}{\Xpn{n}}{K_{n+1}}$, \ where \ 
$K_{n+1}:=\Cok \pi_{n+2}\rho\q{n}$, \ 
up to self-homotopy equivalences of \ $\Xpn{n}$ \ over \ $\BCL$ \ and \ 
$\AutL(K_{n+1})$. \ In particular, this distinguishes the class of \ 
$\pinC{n+1}X$ \ in \ 
$\Ext_{\RM{\Lambda}}(\Ker(Tr\q{n})_{n+1}, \Cok(Tr\q{n})_{n+2})$.
\end{enumerate}
\end{thm}

Note that \ $\Gamma_{n+1}\Xpn{n}=\Gamma_{n+1}\Xpn{n+1}=\Gamma_{n+1}X$, \ 
by Corollary \ref{czero}.

\subsection{Moduli spaces}
\label{sreal}\stepcounter{thm}

It is possible to refine the statement of our fundamental problem of
lifting \ $G\in\D$ \ to $\C$ in terms of \emph{moduli} spaces\vsm : 

Given a model category $\C$, let $\We$ be a homotopically small 
subcategory of $\C$,  such that all maps in $\We$ are weak equivalences,
and if \ $f:X\to Y$ \ is a weak equivalence in $\C$ with either $X$ or
$Y$ in $\We$, then \ $f\in\We$. \  
Recall from \cite[\S 2.1]{DKanCD} that the nerve \ $B\We$ \ of such a
category is called a \emph{classification complex}. Its components are
in one-to-one correspondence with the weak homotopy types (in $\C$) 
of the objects of $\We$, and the component containing \ $X\in\C$ \ is
weakly equivalent to the classifying space \ $B\haut X$ \ of the
monoid of self-weak equivalences of $X$. 

\begin{defn}\label{drealsp}\stepcounter{subsection}
Given a spherical functor \ $T:\C\to\D$ \ and \ $G\in\D$, \ we denote 
by \ $\M(G)$ \ the category of objects in $\D$ weakly equivalent to 
$G$ (with weak equivalences as morphisms), and by \ $\TM(G)$ \ 
the category of objects \ $X\in\C$ \ such that \ 
$TX\in\M(G)$ \ (again, with weak equivalences in $\C$ as 
morphisms). The ``pointed'' version is denoted
by \ $\R{}(G)$ \ -- \ the category of pairs \ $(X,\rho)$, \ 
where \ $X\in\C$ \ and \ $\rho:G\to TX$ \ is a specified weak equivalence.
\end{defn}

In all our examples the obvious functors \ $\R{}(G)\xra{F}\TM(G)\xra{T} \M(G)$ \  
preserve fibrant and cofibrant objects, and thus induce a homotopy 
pullback diagram:
$$
\xymatrix{
B\R{}(G) \rto^{BF} \dto &
	B\TM(G) \dto^{BT} \\
\{\Id_{G}\} \rto & B\M(G)\\
}
$$
\noindent and there are weak equivalences \ 
$B\TM(G) \simeq \coprod_{X\in \pi_{0}\TM(G)} B\haut X$, \ where \ 
$B\M(G) \simeq B\Aut(G)$ \ for \ $\Aut(G)$ \ the monoid 
of self weak equivalences of $G$.

\subsection{Towers of moduli spaces}
\label{streal}\stepcounter{thm}

Although \ $B\TM(G)$ \ is the more natural object of interest in 
our context, it is more convenient to study \ $B\R{}(G)$ \ by means of a 
tower of fibrations, corresponding to the Postnikov system of \ $X\in\R{}(G)$:

Let \ $\R{n}(G)$ \ denote the category whose objects are pairs \ 
$(\Xpn{n},\rho')$, \ where \ $\Xpn{n}\in\C$ \ has \ 
$P_{n+1}\Xpn{n}\simeq\Xpn{n}$ \ and \ 
$\rho':P_{n+1}G\to P_{n+1}T\Xpn{n}$ \ is a weak equivalence.
The maps of \ $\R{n}(G)$ \ are weak equivalences compatible 
with the maps \ $p\q{n}$.	

As in \cite[Thm.\ 9.4]{BDGoeR}, one can show that \ 
$B\R{}(G)\simeq \holim_{n} B\R{n}(G)$, \ so we may try to obtain 
information about the moduli space \ $\TM(G)$ \ by studying the 
successive stages in the tower:
%
\setcounter{equation}{\value{thm}}\stepcounter{subsection}
\begin{equation}\label{esixteen}
\dotsc B\R{n+1}(G) \xra{BF_{n}}B\R{n}(G) \xra{BF_{n-1}}\dotsc
\to B\R{1}(G).
\end{equation}
\setcounter{thm}{\value{equation}}

However, from the discussion above we see that we need several 
intermediate steps in the study of \ $B\R{n+1}(G)\to B\R{n}(G)$, \ 
corresponding to the additional choices made in obtaining \ 
$\bPa{n+1}X$ \ and \ 
$p\q{n+1}:P_{n+2}G\xra{\simeq}P_{n+2}T\bPa{n+1}X$ \ from \ 
$\bPa{n}X$ \ and \ $p\q{n}:P_{n+1}G\xra{\simeq}P_{n+1}T\bPa{n}X$. \ 
As a result one obtains a refinement of the tower \ \eqref{esixteen}, \ 
where the successive fibers $F$ are either empty, or else generalized
Eilenerg-Mac Lane spaces, whose homotopy groups may be described in
terms of appropriate Quillen cohomology groups. We leave the details
to the reader; compare \cite[Thm.\ 9.6]{BDGoeR}.

%
%
\sect{Applying the theory}
\label{cat}

The approach to the lifting problem for a spherical functor \ 
$T:\C\to\D$ \ described in the previous section is somewhat
unwieldy. However, in specific applications it may simplify in various
ways.  We illustrate this by a number of examples:

\subsection{Singular chains}
\label{ssc}\stepcounter{thm}

Consider the singular chain functor \ $C_{\ast}:\Ta\to\Ch$, \ which
in the form \ $T:\G\to s\Abgp$ \ is induced by abelianization 
(see \S \ref{ssfunc}(a)). Thus, given a chain complex \ $G_{\ast}$, \
we would like to find all topological spaces $X$ (if any) with \ 
$C_{\ast}X\simeq G_{\ast}$. \ Over $\bZ$, this is equivalent to the
question of realizing a given sequence of homology groups. 

Our approach uses  Whitehead's exact sequence \
\eqref{etwentyfive} \ to relate the (trivial) Postnikov system for the
chain complex \ $G_{\ast}$ \ to the modified Postnikov system for the
space $X$,  in which we attach at each stage not a single new homotopy
group, but a pair of groups in adjacent dimensions, corresponding to
the image and kernel respectively of the Hurewicz homomorphism.

It should be observed that the functor $T$ involves only 
``algebraic'' categories \ $\C=s\hC$, \ where $\hC$ \ -- \ in our
case, \ $\Gp$ \ or \ $\Abgp$ \ -- \ has a trivial model category
structure, as in \S \ref{ermc}(a-b). \ The analysis 
in Section \ref{cfib} then simplifies considerably, in as much as the 
categories of \PCa s and \PDa s are simply \ $\Gp$ \ and \ $\Abgp$, \
respectively. 

As noted in the Introduction, Baues's \cite[VI, (2.3)]{BauCF} is
actually a generalization the obstruction theory described here for
this case. His earlier approach in \cite{BauHH} (as well as that of 
Benkhalifa in \cite{BenkT} is parallel to this, though not framed in the
same cohomological language. See \cite{MandE} for another viewpoint. 

\subsection{Rationalization}
\label{srat}\stepcounter{thm}

On the other hand, the rationalization functor \
$(-)_{\bQ}:\TT\to\TT_{\bQ}$, \ induced by the completed group ring
functor \ $\hQ:\Gp\to\Hopf$ \ (cf.\ \S \ref{ssfunc}(b)), is spherical
but not special (Def.\ \ref{dspecial}), and so the theory described
here does not apply as is. In fact, one can see why if one
considers the comparison exact sequence for $\hQ$ \ (\S \ref{satt}(b)):
given a (simply-connected) rational space \ $G\in\TT_{\bQ}$, \ for
each $\bQ$-vector space \ $\pi_{n}G$, \ we need an abelian group \  
$A=\pi_{n}X$ \ such that \ $A\otimes\bQ\cong\pi_{n}G$, \ and then lift
the rational $k$-invariants for $X$ to integral ones. \ Thus, much of
the indeterminacy for $X$ is algebraic.

\subsection{Suspension}
\label{ssusp}\stepcounter{thm}

The suspension functor \ $\Sigma:\Ta\to\Ta$, \ induced by the free
group functor \ $\hF:\Seta\to\Gp$ \ as in \S \ref{ssfunc}(c), is
similar to singular chains, with the generalized EHP sequence
replacing the ``certain long exact sequence'', and the modified
Postniov systems involve the kernel and image of the suspension
homomorphism \ $E:\pi_{n}X\to\pi_{n+1}\Sigma X$.

\subsection{Homotopy groups}
\label{shg}\stepcounter{thm}

The motivating example for the treatment in this paper \ -- \ and the
only one which requires the full force of Section \ref{cfib} \ -- \ 
is the functor \ $\pis:\Ta\to\PAlg$, \ prolonged to simplicial spaces
(as in as in \S \ref{ssfunc}(d)). \  However, even this case
simplifies greatly if we want to realize a single \Pa\ $\Lambda$ \ -- \
that is, we take \ $G\in s\PAlg$ \ to be the constant simplicial \Pa \ $\BL$.

Indeed, given a simplicial space $X$ with \ $\pis X\simeq\BL$ \ (which
implies that \ $\pis \|X\|\cong G$), from the spiral exact
sequence \eqref{efifteen} we find that \
$\pin_{n}X\cong\Omega^{n}\Lambda$ \ for all \ $n\geq 0$, \ so that \ 
$h_{n}:\pin_{n}X\to\pin_{n}\pis X$ \ is trivial for \ $n>0$. \ We do
not need the modified Postnikov system in this case: the
obstructions to realizing $\Lambda$ (or $G$) are just the classes \ 
$\chi_{n}\in H^{n+3}(\Lambda;\Omega^{n+1}\Lambda)$, \ and the 
difference obstructions distinguishing between the different
realizations are \ $\delta_{n}\in H^{n+2}(\Lambda;\Omega^{n+1}\Lambda)$ \ 
($n\geq 1$). \ See \cite{BDGoeR} and \cite[\S 5]{BJTurR} for two
descriptions of this case.

\begin{remark}\label{risf}\stepcounter{subsection}
Our obstruction theory is irrelevant, of course, for the inverse
spherical functor \ $U:\G\to\Sa$ \ (see \S \ref{eisf}) \ -- \ that is, 
in determining loop structures on a given topological space.
Nevertheless, from \eqref{etwentythree} we can easily recover the
well-known fact that \ $X\simeq\Omega Y$ \ is a loop space if and only if its
$k$-invariants are suspensions of those of $Y$ (cf.\ \cite{AHKaneN}).
\end{remark}

\subsection{Lifting morphisms}
\label{slm}\stepcounter{thm}

In all of the above examples, one can ask the analogous question
regarding the lifting of \emph{maps}, or more complicated diagrams,
from $\D$ to $\C$. This can be addresses via Theorem \ref{tfour} by
transfering the spherical structure from $\C$ and $\D$ to the diagram
categories \ $\C^{I}$ \ and \ $\D^{I}$ \ (cf.\ \S \ref{ermc}(d)).
See \cite[\S 8]{BJTurR} for a detailed example. 

Note that the $k$-invariants for a map of chain complexes are not
trivial (cf.\ \cite[(3.8)]{DolH}), so the theory for realizing chain maps in \ 
$\Ta$ \ is correspondingly more complicated.


\begin{thebibliography}{DKSt2}
%
\bibitem[AHK]{AHKaneN}
K.~Aoki, E.~Honma, \& T.~Kaneko,
``On natural systems of some spaces'',\hsm
\textit{Proc.\ Jap.\ Acad.} \textbf{32} (1956), pp.~564-567.
%
\bibitem[AC]{ACuH}
M.~Arkowitz \& C.R.~Curjel,
``The Hurewicz homomorphism and finite homotopy invariants'',\hsm
\textit{Trans. AMS} \textbf{110} (1964), pp.~538-551.
%
\bibitem[Ba1]{BauR}
H.-J.~Baues,
``Relationen f\"{u}r prim\"{a}re Homotopieoperationen und eine
  verallgemeinerte $EHP$-Sequenz'',\hsm
\textit{Ann.\ Sc.\ {\'{E}}c.\ Norm.\ Sup.} \textbf{8} (1975), pp.~509-533.
%
\bibitem[Ba2]{BauCHF}
H.-J.~Baues,
\textit{Combinatorial Homotopy and $4$-Dimensional Complexes},\hsm
Gruyter Expositions in Mathematics \textbf{2}, Walter de Gruyter,
Berlin-\-New York, 1991.
%
\bibitem[Ba3]{BauHH}
H.-J.~Baues,
\textit{Homotopy type and homology},\hsm
Oxford University Press, New York, 1996.
%
\bibitem[Ba4]{BauCF}
H.-J.~Baues,
\textit{Combinatorial Foundation of Homology and Homotopy},\hsm
Applications to Spaces, Diagrams, Transformation Groups,
  Compactifications, Differential Algebras, Algebraic Theories, Simplicial
  Objects, and Resolutions.
Springer \textit{Monographs in Mathematics}, Berlin-\-New York, 1999.
%
\bibitem[Be]{BenkT}
M.~Benkhalifa,
``Sur le type d'homotopie d'un {CW}-complexe'',\hsm
\textit{Homology, Homotopy Appl.} \textbf{5} (2003), pp.~101-120.
%
\bibitem[BH]{BHilS}
I.~Berstein \& P.J.~Hilton,
``On suspensions and comultiplications'',\hsm
\textit{Topology} \textbf{2} (1963), pp.~73-82.
%
\bibitem[Bl1]{BlaN}
D.~Blanc,
``New model categories from old'',\hsm
\textit{J.\ Pure \& Appl.\ Alg.} \textbf{109} (1996), pp.~37-60.
%
\bibitem[Bl2]{BlaCW}
D.~Blanc,
``CW simplicial resolutions of spaces, with an application to loop spaces'',\hsm 
\textit{Top.\ \& Appl.} \textbf{100} (2000), pp.~151-175.
%
\bibitem[Bl3]{BlaG}
D.~Blanc,
``Generalized Quillen cohomology'',\hsm 
preprint, 2006.
%
\bibitem[BDG]{BDGoeR}
D.~Blanc, W.G.~Dywer, \& P.G.~Goerss,
``The realization space of a \Pa: a moduli problem in algebraic 
topology'',\hsm \textit{Topology} \textbf{43} (2004), pp.~857-892.
%
\bibitem[BJT]{BJTurR}
D.~Blanc, M.J.~Johnson, \& J.M.~Turner,
``On realizing diagrams of \Pa s'',\hsm
\textit{Algebraic \& Geometric Topology}, to appear.
%
\bibitem[BP]{BPescF}
D.~Blanc \& G.~Peschke,
``The fiber of functors between categories of algebras'',\hsm
\textit{J. Pure Appl. Alg.}, to appear.
%
\bibitem[Bor]{BorcH2}
F.~Borceux,
\textit{Handbook of Categorical Algebra, Vol. 2: Categories and
Structures},\hsm
Encyc. Math. {\&} its Appl., \textbf{51},
Cambridge U. Press, Cambridge, UK, 1994.
%
\bibitem[Bou]{BousCR}
A.K.~Bousfield,
``Cosimplicial resolutions and homotopy spectral sequences in 
model categories'',\hsm
\textit{Geom.\ \& Topology} \textbf{7} (2003), pp.~1001-1053.
%
\bibitem[BK]{BKaH}
A.K.~Bousfield \& D.M.~Kan,
\textit{Homotopy Limits, Completions, and Localizations},\hsm
Springer \textit{Lec.\ Notes Math.} \textbf{304}, Berlin-\-New York, 1972.
%
\bibitem[D]{DolH}
A.~Dold,
``Homology of symmetric products and other functors of complexes'',\hsm
\textit{Ann.\ Math., Ser.\ 2}, \textbf{68} (1958), pp.~54-80.
%
\bibitem[DK1]{DKanCD}
W.G.\ Dwyer \& D.M.\ Kan, 
``A classification theorem for diagrams of simplicial sets'',\hsm
\textit{Topology} \textbf{23} (1984), pp.~139-155.
%
\bibitem[DK2]{DKanO}
W.G.\ Dwyer \& D.M.\ Kan, 
``An obstruction theory for diagrams of simplicial sets'',\hsm
\textit{Proc.\ Kon.\ Ned.\ Akad.\ Wet.\ - Ind.\ Math.} \textbf{46} (1984), 
pp.~139-146.
%
\bibitem[DKS1]{DKStE} 
W.G.~Dwyer, D.M.~Kan, \& C.R.~Stover, 
``An $E^{2}$ model category structure for pointed simplicial spaces'',\hsm 
\textit{J.\ Pure \& Appl.\ Alg.} \textbf{90} (1993), pp.~137-152.
%
\bibitem[DKS2]{DKStB} 
W.G.~Dwyer, D.M.~Kan, \& C.R.~Stover, 
``The bigraded homotopy groups $\pi_{i,j}X$ of a pointed simplicial 
space'',\hsm 
\textit{J.\ Pure Appl.\ Alg.} \textbf{103} (1995), pp.~167-188.
%
\bibitem[G]{GanG}
T.~Ganea,
``A generalization of the homology and homotopy suspensions'',\hsm
\textit{Comm.\ Math.\ Helv.} \textbf{39} (1965), pp.~295-322.
%
\bibitem[GH]{GHopkM}
P.G.~Goerss \& M.J.~Hopkins,
``Moduli spaces of commutative ring spectra'',\hsm 
in A.~Baker and B.~Richter, eds., 
\textit{Structured Ring Spectra}, London Math.\ Soc.\ 
Lect.\ Notes \textbf{315}, Cambridge U.\ Press, Cambridge, 2004, pp.~151-200.
%
\bibitem[Hi]{PHirM}
P.S.~Hirschhorn,
\textit{Model Categories and their Localizations},\hsm
Math.\ Surveys \& Monographs \textbf{99}, AMS, Providence, RI, 2002.
%
\bibitem[Ho]{HopfT}
H.~Hopf,
``\"{U}ber die Topologie der Gruppen-Mannigfaltkeiten und ihre 
Verallgemeinerungen'',\hsm
\textit{Ann.\ Math.\ (2)} \textbf{42} (1941), pp.\ 22-52.
%
\bibitem[J]{JardBE}
J.F.~Jardine,
``Bousfield's $E\sb{2}$ Model Theory for Simplicial Objects'',\hsm
in P.G. Goerss \& S.B. Priddy, eds., 
\textit{Homotopy Theory: Relations with Algebraic Geometry, Group 
Cohomology, and Algebraic {$K$}-Theory},
Contemp.\ Math. \textbf{346}, AMS, Providence, RI 2004, pp.~305-319.
%
\bibitem[Mc]{MacH} 
S.~Mac~Lane, 
\textit{Homology},\hsm 
Springer-Verlag \textit{Grund.\ math.\ Wissens.} \textbf{114}, 
Berlin-\-New York  1963.
%
\bibitem[Man]{MandE}
M.A.~Mandell,
``$E\sb{\infty}$-algebras and $p$-adic homotopy theory'',\hsm
\textit{Topology} \textbf{40} (2001), pp..~43-94.
%
\bibitem[May]{MayS}
J.P.~May,
\textit{Simplicial Objects in Algebraic Topology},\hsm
U.\ Chicago Press, Chicago-\-London, 1967.
%
\bibitem[Ne]{NeiP}
J.A.~Neisendorfer,
\textit{Primary homotopy theory},\hsm
AMS \textit{Memoirs} \textbf{232}, Am. Math. Soc., Providence, RI, 1980.
%
\bibitem[No]{NomE}
Y.~Nomura,
``On extensions of triads'',\hsm
\textit{Nagoya Math.\ J.} \textbf{27} (1966), pp.~249-277
%
\bibitem[Q1]{QuiH}
D.G.~Quillen,
\textit{Homotopical Algebra},\hsm
Springer \textit{Lec.\ Notes Math.} \textbf{43},
Berlin-\-New York, 1963.
%
\bibitem[Q2]{QuiR}
D.G.~Quillen,
``Rational homotopy theory'',\hsm
\textit{Ann.\ Math.} \textbf{90} (1969), pp.~205-295.
%
\bibitem[Q3]{QuiC} 
D.G.~Quillen, 
``On the (co-)homology of commutative rings'',\hsm 
\textit{Applications of Categorical Algebra}, \ Proc.\ Symp.\ Pure Math.\ 
\textbf{17}, AMS, Providence, RI, 1970, pp.~65-87.
%
\bibitem[R]{RobO}
C.A.~Robinson,
``Obstruction theory and the strict associativity of Morava
$K$-theories'',\hsm 
in S.M.~Salamon, B.~Steer, \& W.A. Sutherland, eds., 
\textit{Advances in Homotopy Theory (Cortona, 1988)}, 
London Math.\ Soc.\ Lect.\ Notes \textbf{139}, Cambridge U.\ Press,
Cambridge, 1989, pp.~143-152. 
%
\bibitem[Se]{SerG}
J.-P.~Serre,
``Groupes d'homotopie et classes de groupes ab\'{e}liens'',\hsm
\textit{Ann.\ Math. (2)} \textbf{58} (1953), pp.~258-294.
%
\bibitem[Sm]{JRSmiT1}
J.R.~Smith,
``Topological realizations of chain complexes I: the general theory'',\hsm
\textit{Top.\ \& Appl.} \textbf{22} (1986), pp.~301-313.
%
\bibitem[Sp]{SpalSM}
J.~Spali\'{n}ski,
``Stratified model categories'',\hsm
\textit{Fund.\ Math.} \textbf{178} (2003), pp.~217-236.
%
\bibitem[Sta]{StaH}
J.D.\ Stasheff,
``Homotopy associativity of {$H$}-spaces, I,II'',\hsm
\textit{Trans.\ AMS} \textbf{108} (1963) pp.\ 275-292, 293-312.
%
\bibitem[Ste]{SteCA}
N.E.~Steenrod,
``The cohomology algebra of a space'',\hsm
\textit{Ens. Math.} \textbf{7} (1961), pp.\ 153-178.
%
\bibitem[Sto]{StoV}
C.R.\ Stover,
``A Van Kampen spectral sequence for higher homotopy groups'',\hsm
\textit{Topology} \textbf{29} (1990), pp.~9-26.
%
\bibitem[Sug]{SugG}
M.\ Sugawara,
``A condition that a space is group-like'',\hsm
\textit{Math.\ J.\ Okayama U.} \textbf{7} (1957), pp.~123-149.
%
\bibitem[Sul]{SulG}
D.P.~Sullivan,
\textit{Geometric Topology, Part I. Localization, periodicity, and
  Galois symmetry},\hsm
Massachusetts Institute of Technology, Cambridge, MA, 1970.
%
\bibitem[T]{ThP}
R.~Thom,
``Sur un probl\`{e}me de Steenrod'',\hsm
\textit{C.R. Acad.\ Sci., Paris} \textbf{236} (1953), pp.~1128-1130.
%
\bibitem[W1]{JWhSF}
J.H.C.\ Whitehead,
``On simply connected, {$4$}-dimensional polyhedra'',
\textit{Comm.\ Math.\ Helv.} \textbf{22} (1949), pp.~48-92.
%
\bibitem[W2]{JWhR}
J.H.C.\ Whitehead,
``On the realizability of homotopy groups'',\hsm
\textit{Ann.\ Math.\ (2)} \textbf{50} (1949), pp.~261-263.
%
\bibitem[W3]{JWhSB}
J.H.C.~Whitehead,
``The secondary boundary operator'',\hsm
\textit{Proc.\ Nat.\ Acad.\ Sci.\ USA} \textbf{36} (1950), pp.~55-60.
%
\bibitem[W4]{JWhC}
J.H.C.~Whitehead,
``A certain exact sequence'',\hsm
\textit{Ann.\ Math.\ (2)} \textbf{52} (1950), pp.~51-110.
%
\bibitem[W5]{JWhSH}
J.H.C.~Whitehead,
``Simple homotopy types'',\hsm
\textit{Amer.\ J.\ Math.} \textbf{72} (1952), pp.~1-57, 
%
\end{thebibliography}
\end{document}